\documentclass[10pt]{article}
 
\usepackage{coursE}

\title{Algebra properties for Besov spaces on unimodular Lie groups}

\author{{\Large Joseph {\sc Feneuil}\footnote {The author is  supported by the ANR project ``Harmonic Analysis at its Boundaries'',   ANR-12-BS01-0013-03}} \\
{Institut Fourier, Grenoble} \\
\vspace{1cm}}
\date{\today}

\begin{document}

\maketitle

\begin{abstract}
We consider the Besov space $B^{p,q}_\alpha(G)$ on a unimodular Lie group $G$ equipped with a sublaplacian $\Delta$. 
Using estimates of the heat kernel associated with $\Delta$, we give several characterizations of Besov spaces, and show an algebra property for $B^{p,q}_\alpha(G) \cap L^\infty(G)$ for $\alpha>0$, $1\leq p\leq+\infty$ and $1\leq q\leq +\infty$. 
These results hold for polynomial as well as for exponential volume growth of balls.
\end{abstract}


\tableofcontents

\pagebreak

\section{Introduction and statement of the results}

\begin{minipage}{0.8\linewidth}
 We use the following notations. $A(x) \lesssim B(x)$ means that there exists $C$ independent of $x$ such that $A(x) \leq C \,  B(x)$ for all $x$. $A(x) \simeq B(x)$ means that $A(x) \lesssim B(x)$ and $B(x) \lesssim A(x)$.
The parameters which the constant is independent to will be either obvious from context or recalled. 
\end{minipage}

\subsection{Introduction}

Let $d\in  \N^*$. In $\R^d$, the Besov spaces $B^{p,q}_\alpha(\R^d)$ are obtained by real interpolation of Sobolev spaces and can be defined, for $p,q\in [1,+\infty]$ and $\alpha\in \R$, 
as the subset of distributions $\mathcal S'(\R^d)$ satisfying
\begin{equation} \label{Fen3Bpqalpha} \|f\|_{B^{p,q}_\alpha} : = \|\psi * f\|_{L^p} + \left(\sum_{k=1}^\infty \left[2^{k\alpha} \|\varphi_k * f\|_{L^p} \right]^q \right)^\frac1q < +\infty\end{equation}
where, if $\varphi \in \mathcal S(\R^d)$ is supported in $B(0,2)\backslash B(0,\frac12)$, $\varphi_k$ and $\psi$ are such that $\mathcal F \varphi_k(\xi) = \varphi(2^{-k}\xi)$ and
$\mathcal F \psi (\xi) = 1 - \sum_{k=1}^\infty \varphi (2^{-k}\xi)$. 

 The norm of the Besov space $B^{p,q}_\alpha(\R^d)$ can be also written by using the heat operator. Indeed, Triebel proved in \cite[Section 2.12.2]{triebelheat,Tri2} that for all $p,q\in [1,+\infty]$, all $\alpha>0$ and all integer $m> \frac\alpha 2$,
\begin{equation} \label{Fen3Bpqalphawithht} \|f\|_{B^{p,q}_\alpha} \simeq \|f\|_{L^p} + \left(\int_{0}^{\infty} t^{(m-\frac{\alpha}{2})q} \left\|\frac{\dr^M H_t}{\dr t^M} f\right\|_{L^p}^q \frac{dt}{t} \right)^\frac1q \end{equation}
where $H_t = e^{t\Delta}$ is the heat semigroup (generated by $-\Delta$). 
Note that we can give a similar characterization by using, instead of the heat semigroup, the harmonic extension or another extensions obtained by convolution (see \cite{uspenskii,mr}).

Another characterization in term of functional using differences of functions was done. 
Define for $M \in \N^*$, $f\in L^p(\R^d)$, $x,h\in \R^d$ the term
$$\nabla^M_h f(x) = \sum_{l=0}^M \left(\begin{matrix} M \\ l \end{matrix} \right) (-1)^{M-l} f(x+lh)$$
and then for $M>\alpha>0$, $p,q\in [1,+\infty]$
\begin{equation} \label{Fen3Strichartztype}
 S_{\alpha,M}^{p,q} f = \left( \int_{\R^d} |h|^{-\alpha q} \|\nabla_h^M f\|_{L^p}^q \right)^\frac1q.
\end{equation}
We have then for all $\alpha>0$, $p,q\in [1,+\infty]$ and $M\in \N$ with $M>\alpha$, 
\begin{equation} \label{Fen32ndcharacterization} \|f\|_{B^{p,q}_\alpha} \simeq \|f\|_{L^p} + S_{\alpha,M}^{p,q} f.\end{equation}

\smallskip

One of the remarkable property of Besov spaces (see \cite[Proposition 1.4.3]{danchin2012fourier}, \cite[Theorem 2, p. 336]{runstalgebra}, \cite[Proposition 6.2]{mr}) is that $B^{p,q}_\alpha(\R^d)\cap L^\infty(\R^d)$ is an algebra for the pointwise product, 
that is for all $\alpha >0$, all $p,q\in [1,+\infty]$, one has
\begin{equation} \label{Fen3AlgebraProperty}
\|fg\|_{B^{p,q}_\alpha} \lesssim \|f\|_{B^{p,q}_\alpha} \|g\|_{L\infty} + \|f\|_{L^\infty} \|g\|_{B^{p,q}_\alpha}.
\end{equation}
The idea of \cite{danchin2012fourier} consists in decomposing the product $fg$ by some paraproducts.
The authors of \cite{mr} wrote $B^{p,q}_\alpha(\R^d)$ as a trace of some weighted (non fractional) Sobolev spaces, 
and thus deduced the algebra property $B^{p,q}_\alpha(\R^d) \cap L^\infty(\R^d)$ from the one of $W^{p,k}(\R^d) \cap L^\infty(\R^d)$.
Notice also that, when $\alpha \in (0,1)$ and $M=1$, the algebra property of $B^{p,q}_\alpha(\R^d) \cap L^\infty(\R^d)$ is a simple consequence of \eqref{Fen32ndcharacterization}.

\smallskip

The property \eqref{Fen3AlgebraProperty} have also been studied in the more general setting of Besov spaces on Lie groups.
Gallagher and Sire stated  in \cite{GS} an algebra property for Besov spaces on $H$-type groups, which are a subclass of Carnot groups. 
In order to do this, they used a some paradifferential calculus and a Fourier transform adapted to $H$-groups.

Moreover, in the more general case where $G$ is a unimodular Lie group with polynomial growth,
they used the definition of Besov spaces obtained using Littlewood-Paley decomposition proved in \cite{FMV}. 
When $\alpha \in (0,1)$, they proved a equivalence of the Besov norms with some functionals using differences of functions, 
in the spirit of \eqref{Fen3Strichartztype}, and thus they obtained an algebra property for $B^{p,q}_\alpha(G) \cap L^\infty(G)$. 
They shows a recursive definition of Besov spaces and wanted to use it to extend the property \eqref{Fen3AlgebraProperty} to $\alpha\geq 1$.
However, it seems to us that there is a small gap in their proof and they actually proved the property  
$\|fg\|_{B^{p,q}_\alpha} \lesssim (\|f\|_{B^{p,q}_\alpha} +\|f\|_{L^\infty})(\|g\|_{B^{p,q}_\alpha} +\|g\|_{L^\infty})$.

\smallskip

In our paper, we defined Besov spaces on unimodular Lie group (that can be of exponential growth) for all $\alpha>0$, and then we proved an algebra property on them.
We used two approaches. One with functionals in the spirit of \eqref{Fen3Strichartztype} and the other one using paraproducts.
We did not state any results on homogeneous Besov spaces because the definition of these spaces need a particular work (a extension of the work in \cite{BDY} to $\alpha\notin (-1,1)$ should work). 
However, we have no doubt that our methods work once we get the proper definition of homogeneous Besov spaces with some good Calder\'on-Zygmund formulas.

\smallskip

Note that methods used in \cite{GS} or in the present paper are similar to the ones in \cite{CRTN} and \cite{BBR}, where fractional Sobolev spaces $L^p_\alpha(G)$ are considered on unimodular Lie groups (and on Riemannian manifolds). 
In these two last articles, the authors proved the algebra property for $L^p_\alpha(G) \cap L^\infty(G)$ when $p\in (1,+\infty)$ and $\alpha>0$.

\subsection{Lie group structure}

In this paper, $G$ is a unimodular connected Lie group endowed with its Haar measure $dx$. We recall that ``unimodular'' means that $dx$ is both left- and right-invariant. 
We denote by $\mathcal L$ the Lie algebra of $G$ and we consider a family $\mathbb X = \{X_1,\dots, X_k\}$ of left-invariant vector fields on $G$ satisfying the Hörmander condition 
(which means that the Lie algebra generated by the family $\mathbb X$ is $\mathcal L$). Denote $\ds \mathcal I_\infty(\N) = \bigcup_{l\in \N} \{1,\dots,k\}^l$.
Then if $I = (i_1,\dots,i_n) \in \mathcal I_\infty(\N)$, the length of $I$ will be denoted by $|I|$ and is equal to $n$, whereas $X_I$ denotes the vector field $X_{i_1}\dots X_{i_n}$.

A standard metric, called the Carnot-Caratheodory metric, is naturally associated with $(G,\mathbb X)$ and is defined as follows. Let $l:\, [0,1]\to G$ be an absolutely continuous path. 
We say that $l$ is admissible if there exist measurable functions $a_1,\dots, a_k: \, [0,1] \to \C$ such that
$$l'(t) = \sum_{i=1}^k a_i(t) X_i(l(t)) \qquad \text{ for a.e. } t\in [0,1].$$
If $l$ is admissible, its length is defined by $\ds |l| = \int_0^1 \left(\sum_{i=0}^k |a_i(t)|^2 \right)^\frac{1}{2} dt$. 
For any $x,y\in G$, the distance $d(x,y)$ between $x$ and $y$ is then the infimum of the lengths of all admissible curves joining $x$ to $y$ (such a curve exists thanks to the Hörmander condition).
The left-invariance of the $X_i$'s implies the left-invariance of $d$. For short, $|x|$ denotes the distance between the neutral $e$ and $x$, and therefore $d(x,y) = |y^{-1}x|$ for all $x$ and $y$ in $G$.

For $r>0$ and $x\in G$, we denote by $B(x,r)$ the open ball with respect to the Carnot-Caratheodory metric centered at $x$ and of radius $r$. 
Define also by $V(r)$ the Haar measure of any ball of radius $r$.

From now and abusively, we will write $G$ for $(G,\mathbb X,d,dx)$. Recall that $G$ has a local dimension (see \cite{NSW}):

\begin{prop} \label{Fen3localdimension}
 Let $G$ be a unimodular Lie group and $\mathbb X$ be a family of left-invariant vector fields satisfying the H\"ormander condition.
Then $G$ has the local doubling property, that is there exists $C>0$ such that 
$$V(2r)\leq CV(r) \qquad \forall 0<r\leq1.$$
More precisely, there exist $d \in \N^*$ and $c,C>0$ such that
$$ cr^d \leq V(r) \leq Cr^d \qquad \forall 0<r\leq 1.$$
\end{prop}

For balls with radius bigger than 1, we have the result of Guivarc'h (see \cite{Guivarch}):

\begin{prop} \label{Fen3infinitedimension}
 If $G$ is a unimodular Lie group, only two situations may occur.

Either $G$ has polynomial growth and there exist $D\in \N^*$ and $c,C>0$ such that
$$ cr^D \leq V(r) \leq Cr^D \qquad \forall r\geq 1,$$

or $G$ has exponential growth and there exist $c_1,c_2,C_1,C_2>0$ such that 
$$c_1e^{c_2r} \leq V(r) \leq C_1 e^{C_2r} \qquad \forall r\geq 1.$$
\end{prop}

We consider the positive sublaplacian $\Delta$ on $G$ defined by
$$\Delta = - \sum_{i=1}^k X_i^2.$$
We will denote by  $H_t = e^{-t\Delta}$ the heat semigroup on $G$ associated with $\Delta$.

\subsection{Definition of Besov spaces}

\begin{defi}
Let $G$ be a unimodular Lie group. We define the Schwartz space $\mathcal S(G)$ as the space of functions $\varphi  \in C^\infty(G)$ where all the seminorms
$$N_{I,c}(\varphi) = \sup_{x\in G} e^{c|x|} |X_I \varphi(x)| \qquad c\in \N, \, I \in \mathcal I_\infty(\N)$$
are finite.
 
The space $\mathcal S'(G)$ is defined as the dual space of $S(G)$.
\end{defi}

\begin{rmq}
Note that we have the inclusion $\mathcal S(G) \subset L^p(G)$ for any $p\in [1,+\infty]$. As a consequence, $L^p(G) \subset \mathcal S'(G)$.  
\end{rmq}

\begin{defi}
Let $G$ be a unimodular Lie group and let $\alpha \geq 0$, $p,q\in [1,+\infty]$. 
The space $f \in B_\alpha^{p,q}(G)$ is defined as the subspace of $\mathcal S'(G)$ made of distributions $f$ such that, for all $t\in (0,1)$, $\Delta^mH_tf\in L^p(G)$ and satisfying
$$\|f\|_{B_\alpha^{p,q}} : = \Lambda^{p,q}_\alpha f + \|H_\frac12 f\|_p < +\infty,$$ 
where 
$$\Lambda^{p,q}_\alpha f : = \left( \int_0^1 \left( t^{m-\frac{\alpha}{2}} \left\| \Delta^m H_t f\right\|_p\right)^q \frac{dt}{t}\right)^{\frac{1}{q}}$$
if $q<+\infty$ (with the usual modification if $q=+\infty$) and $m$ stands for the only integer such that $\frac{\alpha}{2} < m \leq \frac{\alpha}{2}+1$.
\end{defi}

\begin{rmq}
Lemma \ref{Fen3htisschwartz} provides that the heat kernel $h_t$ is in $\mathcal S(G)$ for all $t>0$. Thus $H_t \varphi \in \mathcal S(G)$ whenever $t>0$ and $\varphi \in \mathcal S(G)$.
When $f \in \mathcal S'(G)$, the term $X_I H_t f$ denotes the distribution in $\mathcal S'(G)$ defined by
$$\left<X_I H_t f,\varphi \right> = (-1)^{|I|}\left<f, H_tX_I\varphi \right> \qquad \forall \varphi \in \mathcal S(G).$$
\end{rmq}

\subsection{Statement of the results}

\begin{prop} \label{Fen3AnalyticityofHt}
Let $G$ be a unimodular Lie group. The one has for all $p\in [1,+\infty]$, all multi indexes $I \in \mathcal I_\infty(\N)$ and all $t\in (0,1)$,
$$\|X_I H_t f\|_p \leq C_{I} t^{-\frac{|I|}{2}} \|f\|_p \qquad \forall f\in L^p(G) \, .$$
\end{prop}

\begin{rmq}
In particular, one has that $\|t\Delta H_t\|_p \lesssim 1$ once $t\in (0,1)$ and for all $p\in [1,+\infty]$. 
When $p\in (1,+\infty)$, since $\Delta$ is analytic on $L^2$ (and thus on $L^p$), 
we actually have $\|t\Delta H_t\|_p \lesssim 1$ for all $t>0$.
The case  
\end{rmq}

The following result gives equivalent definitions of the Besov spaces $B^{p,q}_\alpha$ only involving the Laplacian.

\begin{theo} \label{FEN3MAIN1}
Let $G$ be a unimodular Lie group and $p,q\in [1,+\infty]$ and $\alpha\geq 0$. 

If $m>\frac\alpha2$ and $t_0$ a real in $\left\{\begin{array}{ll} (0,1) & \text{if } \alpha=0 \\ \left[0,1\right) & \text{if } \alpha>0 \end{array}\right.$,
then the following norms are equivalent to the norm of $B^{p,q}_\alpha(G)$. 
\begin{enumerate}[(i)]
\item $\ds \left( \int_0^1 \left( t^{m-\frac{\alpha}{2}} \left\| \Delta^m H_t f\right\|_p\right)^q \frac{dt}{t}\right)^{\frac{1}{q}} + \|H_{t_0}f\|_p$.

\item $\ds\|H_{t_0}f\|_p +  \left(\sum_{j\leq -1} \left[2^{j(m-\frac{\alpha}{2})} \|\Delta^m H_{2^{j}} f\|_p \right]^q \right)^{\frac1q}$.

\item $\ds \|H_{t_0}f\|_p + \left(\sum_{j\leq -1} \left[2^{-j\frac{\alpha}{2}} \left\|\int_{2^j}^{2^{j+1}}\left|(t\Delta)^m H_{t} f\right| \frac{dt}{t}\right\|_p \right]^q \right)^{\frac1q}$

if we assume that $\alpha>0$.
\end{enumerate}
\end{theo}

\begin{rmq}
Here and after, we say that ``a norm $N$ is equivalent to the norm in $B^{p,q}_\alpha$'' if and only if the space of distributions $f\in \mathcal S^{\prime}$ such that $\Delta^mH_tf$ is a locally integrable function in $G$ for all $t>0$ and $N(f)<+\infty$ coincides with $B^{p,q}_\alpha$ and the norm $N$ is equivalent to $\left\Vert \cdot \right\Vert_{B^{p,q}_\alpha}$.
\end{rmq}

The previous theorem allows us to recover some well known facts about Besov spaces in $\R^d$.

\begin{cor} \label{Fen3InclusionofBpqs} [Embeddings]
Let $G$ a unimodular Lie group, $p,q,r\in [1,+\infty]$ and $\alpha\geq 0$. 
We have the following continuous embedding
$$B^{p,q}_\alpha(G) \subset B^{p,r}_\alpha(G)$$
once $q\leq r$.
\end{cor}

\begin{cor} \label{Fen3InterpolationTh} [Interpolation]

Let $G$ be a unimodular Lie group. 
Let $s_0,s_1\geq 0$ and $1\leq p_0,p_1,q_0,q_1\leq \infty$.

\noindent Define
$$s^* = (1-\theta)s_0 + \theta s_1$$
$$\frac{1}{p^*} = \frac{1-\theta}{p_0} + \frac{\theta}{p_1}$$
$$\frac{1}{q^*} = \frac{1-\theta}{q_0} + \frac{\theta}{q_1}.$$

The Besov spaces form a scale of interpolation for the complex method, that is, if $s_0 \neq s_1$,
$$(B^{p_0,q_0}_{s_0}, B^{p_1,q_1}_{s_1})_{[\theta]} = B^{p^*,q^*}_{s^*}.$$
\end{cor}

The following result is another characterization of Besov spaces, using explicitly the family of vector fields $\mathbb X$.

\begin{theo}
\label{Fen3DiscreteEquivalenceD}
Let $G$ be a unimodular Lie group, $p,q\in [1,+\infty]$ and $\alpha > 0$. Let $\bar m$ be an integer strictly greater than $\alpha$. Then
\begin{equation} \label{Fen3discreteXsup} \|H_\frac12f\|_p + \left(\sum_{j\leq -1} \left[2^{j\frac{\bar m - \alpha}{2}} \max_{t\in [2^j,2^{j+1}]}\sup_{|I|\leq \bar m}\|X_I H_{t} f\|_p \right]^q \right)^{\frac1q}\end{equation}
is an equivalent norm in $B^{p,q}_\alpha(G)$.
\end{theo}

With the use of paraproducts, we can deduce from Corollary \ref{Fen3InterpolationTh} and Theorem \ref{Fen3DiscreteEquivalenceD} the complete following Leibniz rule.

\begin{theo} \label{FEN3MAIN2}
 Let $G$ be a unimodular Lie group, $0<\alpha$ and $p,p_1,p_2,p_3,p_4,q\in [1,+\infty]$ such that
$$\frac{1}{p_1} + \frac{1}{p_2} = \frac{1}{p_3} + \frac{1}{p_4} = \frac{1}{p}.$$ 
Then for all $f\in B^{p_1,q}_\alpha \cap L^{p_3}$ and all $g \in B^{p_4,q}_\alpha \cap L^{p_2}$, one has
\begin{equation}\label{Fen3LeibnizRule}\|fg\|_{B^{p,q}_\alpha} \lesssim \|f\|_{B^{p_1,q}_\alpha}\|g\|_{L^{p_2}} + \|f\|_{L^{p_3}} \|g\|_{B^{p_4,q}_\alpha}.\end{equation}
\end{theo}

\begin{rmq}
The Leibniz rule implies that $B^{p,q}_\alpha(G) \cap L^\infty(G)$ is an algebra under pointwise product, that is
$$\|fg\|_{B^{p,q}_\alpha} \lesssim \|f\|_{B^{p,q}_\alpha}\|g\|_{L^{\infty}} + \|f\|_{L^{\infty}} \|g\|_{B^{p,q}_\alpha}.$$
\end{rmq}

Let us state another characterization of $B^{p,q}_\alpha$ in term of functionals using differences of functions.

Define $\nabla_y f(x) = f(xy)-f(x)$ for all functions $f$ on $G$ and all $x,y\in G$. 
Consider the following sublinear functional
$$L^{p,q}_\alpha(f) = \left(\int_{|y|\leq 1} \left(\frac{\|\nabla_y f\|_p}{|y|^\alpha}\right) ^q \frac{dy}{V(|y|)}\right)^{\frac{1}{q}}.$$

\begin{theo} \label{Fen3Main3}
Let $G$ be a unimodular Lie group. Let $p,q\in [1,+\infty]$. Then for all $f\in L^p(G)$,
$$L^{p,q}_\alpha(f) + \|f\|_p \simeq \Lambda^{p,q}_\alpha(f) + \|f\|_p $$
once $\alpha \in (0,1)$. 
\end{theo}

\begin{rmq}
When $G$ has polynomial volume growth, Theorem \ref{Fen3Main3} is the inhomogeneous counterpart of Theorem 2 in \cite{SC2}. 
Note that this statement is new when $G$ has exponential volume growth.
\end{rmq}

\begin{rmq}
From Theorem \ref{Fen3Main3}, we can deduce the Leibniz rule stated in Theorem \ref{FEN3MAIN2} in the case $\alpha \in (0,1)$. 
\end{rmq}

As Sobolev spaces, Besov spaces can be characterized recursively.

\begin{theo} \label{Fen3Main4}
 Let $G$ be a unimodular Lie group. Let $p,q\in [1,+\infty]$ and $\alpha>0$. Then
$$f\in B^{p,q}_{\alpha+1}(G) \Leftrightarrow \forall i, \, X_if \in B^{p,q}_\alpha(G) \text{ and } f\in L^p(G).$$
\end{theo}

\begin{rmq}
Note that a similar statement is established in \cite{GS}. 
However, we prove this fact for $p\in [1,+\infty]$ while the authors of \cite{GS} used the boundedness of the Riesz transforms and thus are restricted to $p \in (1,+\infty)$.
\end{rmq}

\section{Estimates of the heat semigroup}

\subsection{Preliminaries}
The following lemma is easily checked:
\begin{lem} \label{Fen3convolutionlemma}
 Let $(A,dx)$ and $(B,dy)$ be two measured spaces. Let $K(x,y): A\times B \to \R_+$ be such that
$$\sup_{x\in A} \int_B K(x,y)dy \leq C_B$$
and
$$\sup_{y\in B} \int_A K(x,y)dx \leq C_A.$$
Let $q\in [1,+\infty]$. Then for all $f\in L^q(B)$
$$\left( \int_A \left\vert \int_B K(x,y)f(y)dy \right\vert^q dx \right)^\frac1q \leq C_B^{1-\frac{1}{q}} C_A^\frac{1}{q} \|f\|_q,$$
with obvious modifications when $q = +\infty$.
\end{lem}

\begin{lem} \label{Fen3Bcalculus}
Let $(a,b)\in (\Z\cup\{\pm \infty\})^2$ such that $a<b$, $0<\alpha<\beta$ two real numbers and $q\in [1,+\infty]$. 
Then there exists $C_{\alpha,\beta}>0$ such that for any sequence $(c_n)_{n\in \Z}$, one has
$$\sum_{j=a}^b \left[ 2^{j\alpha} \sum_{n=a}^b 2^{-\max\{n,j\}\beta} c_n \right]^q \lesssim \sum_{n=a}^b \left[2^{(\alpha-\beta)n} c_n\right]^q.$$
\end{lem}

\begin{dem}
We have 
$$\sum_{j=a}^b \left[ 2^{j\alpha} \sum_{n=a}^b 2^{-\max\{n,j\}\beta} c_n \right]^q = \sum_{j=a}^b \left[ \sum_{n=a}^b K(n,j) d_n \right]^q$$
with $d_n = 2^{n(\alpha-\beta)} c_n$ and $K(n,j) = 2^{(j-n)\alpha} 2^{(n-\max\{j,n\})\beta}$. 

According to Lemma \ref{Fen3convolutionlemma}, one has to check that
$$\sup_{j\in \bb a,b\bn} \sum_{n=a}^b K(n,j)\lesssim 1$$
and
$$\sup_{n\in \bb a,b\bn} \sum_{j=a}^b K(n,j) \lesssim 1.$$
For the first estimate, check that
\[\begin{split}
   \sup_{j\in \bb a,b\bn} \sum_{n=a}^b K(n,j) & = \sup_{j\in \bb a,b\bn} \left[ 2^{j(\alpha-\beta)} \sum_{n=a}^j 2^{n(\beta-\alpha)} + 2^{j\alpha} \sum_{n=j+1}^{b} 2^{-n\alpha} \right] \\
& \leq \sup_{j\in\Z} \left[ 2^{j(\alpha-\beta)} \sum_{n=-\infty}^j 2^{n(\beta-\alpha)} + 2^{j\alpha} \sum_{n=j+1}^{+\infty} 2^{-n\alpha} \right] \\
& \lesssim 1,
  \end{split}\]
since $\beta-\alpha>0$ and $\alpha>0$.

The second estimate can be checked similarly:
\[\begin{split}
   \sup_{n\in \bb a,b\bn} \sum_{j=a}^b K(n,j) & = \sup_{j\in \bb a,b\bn} \left[ 2^{-n\alpha} \sum_{j=a}^n 2^{j\alpha} + 2^{n(\beta-\alpha)} \sum_{j=n+1}^{b} 2^{j(\alpha-\beta)} \right] \\
& \lesssim 1.
  \end{split}\]
\end{dem}

\begin{prop} \label{Fen3convolutionbyht}
Let $s\geq 0$ and $c>0$. Define, for all $t\in (0,1)$ and all $x,y\in G$, 
$$K_t(x,y) = \left( \frac{|y^{-1}x|^2}{t} \right)^s \frac{1}{V(\sqrt{t})} e^{-c\frac{|y^{-1}x|^2}{t}}.$$
Then, for all $q\in [1,+\infty]$,
$$\left(\int_G \left( \int_G K_t(x,y) g(y) dy \right)^q dx \right)^\frac{1}{q} \lesssim \|g\|_q.$$ 
\end{prop}

\begin{dem}
Let us check that the assumptions of Lemma \ref{Fen3convolutionlemma} are satisfied. For all $x\in G$ and all $t\in (0,1)$, 
\[\begin{split}
   \int_G K_t(x,y) dy & = \frac{1}{V(\sqrt{t})} \int_{G} \left(\frac{|y^{-1}x|^2}{t} \right)^s e^{-c\frac{|y^{-1}x|^2}{t}} dy \\
& = \frac{1}{V(\sqrt{t})} \int_{|y^{-1}x|^2<t} \left(\frac{|y^{-1}x|^2}{t} \right)^s e^{-c\frac{|y^{-1}x|^2}{t}} dy \\
& \qquad + \frac{1}{V(\sqrt{t})} \int_{|y^{-1}x|^2\geq t} \left(\frac{|y^{-1}x|^2}{t} \right)^s e^{-c\frac{|y^{-1}x|^2}{t}} dy \\
& = I_1 +I_2.
\end{split}\]
The term $I_1$ is easily dominated by 1. As for $I_2$, it is estimated as follows:
\[\begin{split}
   I_2 & = \sum_{j=0}^\infty \frac{1}{V(\sqrt{t})} \int_{2^j \sqrt{t} \leq |y^{-1}x| < 2^{j+1}\sqrt{t}} \left(\frac{|y^{-1}x|^2}{t} \right)^s e^{-c\frac{|y^{-1}x|^2}{t}} dy \\
& \lesssim \sum_{j=0}^\infty \frac{V(2^{j+1}\sqrt{t})}{V(\sqrt{t})} 4^{js} e^{-c4^j}.
  \end{split}\]
Notice that Propositions \ref{Fen3localdimension} and \ref{Fen3infinitedimension} imply that $\frac{V(2^{j+1}\sqrt{t})}{V(\sqrt{t})} \lesssim 2^{jd}$ if $2^{j}\sqrt{t} \leq 1$ and 
\begin{equation} \label{Fen3croissanceexp}\frac{V(2^{j+1}\sqrt{t})}{V(\sqrt{t})} = \frac{V(2^{j+1}\sqrt{t})}{V(1)} \frac{V(1)}{V(\sqrt{t})} \lesssim e^{C2^j} 2^{jd}\end{equation}
if $2^{j}\sqrt{t} \geq 1$. Hence,
$$\sum_{j = 0}^\infty \frac{V(2^{j+1}\sqrt{t})}{V(\sqrt{t})} 4^{js} e^{-c4^j} \lesssim \sum_{j=0}^\infty e^{-c' 4^j} \lesssim 1,$$
which yields with the uniform estimate
\begin{equation} \label{Fen3KisL1} \int_G K_t(x,y) dy \lesssim 1.\end{equation}

In the same way, one has
$$\int_G K_t(x,y) dx \lesssim 1.$$

Lemma \ref{Fen3convolutionlemma} provides then the desired result.
\end{dem}

\begin{prop} \label{Fen3convolutionbyhtbis}
Let $s\geq 0$ and $c>0$. Define
$$K(t,y) = \left( \frac{|y|^2}{t} \right)^s \frac{V(|y|)}{V(\sqrt{t})} e^{-c\frac{|y|^2}{t}}.$$
Then, for all $q\in [1,+\infty]$,
$$\left(\int_0^1 \left( \int_G K(t,y) g(y) dy \right)^q \frac{dt}{t} \right)^\frac{1}{q} \lesssim \left( \int_G |g(y)| \frac{dy}{V(|y|)} \right)^\frac1q.$$ 
\end{prop}

\begin{dem}
 Let us check again that the assumptions of Lemma \ref{Fen3convolutionlemma} are satisfied, that are in our case
$$\sup_{t\in (0,1) } \int_G K(t,y)dy \leq C_B$$
and
$$\sup_{y\in G} \int_0^1 K(t,y) \frac{dt}{t} \leq C_A.$$

The first one is exactly as the estimate \eqref{Fen3KisL1}.For the second one, check that
\[\begin{split}
   \int_0^1 K(t,y) \frac{dt}{t} & = \int_{0}^1  \frac{V(|y|)}{V(\sqrt{t})} \left(\frac{|y|^2}{t} \right)^s e^{-c'\frac{|y|^2}{t}} \frac{dt}{t} \\
& \lesssim \int_{0}^{|y|^2} \frac{V(|y|)}{V(\sqrt{t})}\left(\frac{|y|^2}{t} \right)^s e^{-c'\frac{|y|^2}{t}} \frac{dt}{t}  + \int_{|y|^2}^\infty \left(\frac{|y|^2}{t}\right)^s \frac{dt}{t}\\
& \lesssim \sum_{j=0}^{+\infty} \int_{4^{-(j+1)}|y|^2}^{4^{-j}|y|^2} \frac{V(|y|)}{V(\sqrt{t})}\left(\frac{|y|^2}{t} \right)^s e^{-c'\frac{|y|^2}{t}} \frac{dt}{t} + 1 \\
& \lesssim \sum_{j=0}^{+\infty} \frac{V(|y|)}{V(2^{-j+1}|y|)} 4^{js}  e^{-c4^j} + 1  \\
& \lesssim \sum_{j=0}^{+\infty} 2^{j(d+2s)} e^{C2^j}  e^{-c4^j} + 1  \\
& \lesssim 1,
\end{split}\]
where the last but one line is obtained with the estimate \eqref{Fen3croissanceexp}. 
\end{dem}

\subsection{Estimates for the semigroup}

Because of left-invariance of $\Delta$ and hypoellipticity of $\frac{\partial }{\partial t} + \Delta$, $H_t = e^{-t\Delta}$ has a convolution kernel $h_t \in C^\infty(G)$ satisfying, for all $f\in L^1(G)$ and all $x\in G$,
$$H_t f(x) = \int_G h_t(y^{-1}x) f(y) dy  = \int_G h_t(y) f(xy) dy = \int_G h_t(y) f(xy^{-1}) dy.$$

The kernel $h_t$ satisfies the following pointwise estimates.

\begin{prop} \label{Fen3htestimates}
Let $G$ be a unimodular Lie group. For all $I \in \mathcal I_\infty(\N)$, there exist$C_I,c_I>0$ such that for all $x\in G$, all $t\in (0,1]$, one has
$$ |X_I h_t(x)| \leq \frac{C_I}{t^{\frac{|I|}{2}}V(\sqrt{t})} \exp\left(-c_I\frac{|x|^2}{t} \right).$$
\end{prop}

\begin{dem}
It is a straightforward consequence of Theorems VIII.2.4, VIII.4.3 and V.4.2. in \cite{VSCC}.
\end{dem}

\begin{lem} \label{Fen3htisschwartz}
 Let $G$ be a unimodular group. Then $h_t \in \mathcal S(G)$ for all $t>0$.
\end{lem}

\begin{dem}
The case $t<1$ is a consequence of the estimates on $h_t$. For $t\geq 1$, just notice that $\mathcal S(G) * \mathcal S(G) \subset \mathcal S(G)$. 
\end{dem}

\begin{prop} \label{Fen3AnalycityofHt}
For all $I \in \mathcal I_\infty(\N)$ and all $p\in [1,+\infty]$, one has
$$\|X_I H_t f\|_p \lesssim t^{-\frac{|I|}{2}} \|f\|_p \qquad \forall t\in (0,1], \, \forall f\in L^p(G).$$
\end{prop}

\begin{dem}
Proposition \ref{Fen3htestimates} yields for any $t\in(0,1]$
\[\begin{split}
   \|X_I H_t f\|_p & \lesssim t^{-\frac{|I|}{2}} \left(\int_G \left|\int_G K_t(x,y) f(y) dy\right|^p dx \right)^\frac1p
  \end{split}\]
where $K_t(x,y) = \frac{1}{V(\sqrt{t})}\exp\left(-c \frac{|y^{-1}x|^2}{t}\right)$.

The conclusion of Proposition \ref{Fen3AnalycityofHt} is an immediate consequence of Proposition \ref{Fen3convolutionbyht}.
\end{dem}

\section{Littlewood-Paley decomposition}

We need a Littlewood-Paley decomposition adapted to this context. In \cite{GS}, the authors used the Littlewood-Paley decomposition proven in \cite[Proposition 4.1]{FMV}, only established in the case of polynomial volume growth. We state here a slightly different version of the Littlewood-Paley decomposition, also valid for the case of exponential volume growth.

\begin{lem} \label{Fen3IdentityFormula}
Let $G$ be a unimodular group and let $m\in \N^*$. For any $\varphi \in \mathcal S(G)$ and any $f\in \mathcal S'(G)$, one has the identities
$$\varphi = \frac{1}{(m-1)!} \int_0^{1} (t\Delta)^{m} H_t \varphi \frac{dt}{t} +  \sum_{k=0}^{m-1} \frac{1}{k!} \Delta^k H_1 \varphi$$
where the integral converges in $\mathcal S(G)$, and
$$f = \frac{1}{(m-1)!} \int_0^{1} (t\Delta)^{m} H_t f \frac{dt}{t} +  \sum_{k=0}^{m-1} \frac{1}{k!} \Delta^k H_1 f$$
where the integral converges in $\mathcal S'(G)$.
\end{lem}

\begin{dem}
We only have to prove the first identity since the second one can be obtained by duality.

Let $\varphi \in \mathcal S(G)$. Check first the formula
$$(m-1)! = \int_0^{+\infty} (tu)^{m} e^{-tu} \frac{dt}{t} = \int_0^1 (tu)^{m} e^{-tu} \frac{dt}{t} + \sum_{k=0}^{m-1} \frac{(m-1)!}{k!} u^k e^{-u}.$$
Thus by functional calculus, since $\varphi \subset L^2(G)$, one has
\begin{equation}\label{Fen3identityformula}\varphi = \frac{1}{(m-1)!} \int_0^1 (t\Delta)^{m} H_t \varphi \frac{dt}{t} + \sum_{k=0}^{m-1} \frac{1}{k!} \Delta^k H_1 \varphi,\end{equation}
where the integral converges in $L^2(G)$. 
Since the kernel $h_t$ of $H_t$ is in $\mathcal S(G)$ for any $t>0$ (see Lemma \ref{Fen3htisschwartz}), 
the formula \eqref{Fen3identityformula} will be proven if we have for any $c\in \N$ and any $I \in \mathcal I_\infty(\N)$,
\begin{equation} \label{Fen3schwartzconvergence}
 \lim_{u\to 0} N_{I,c}\left( \int_0^u (t\Delta)^{m} H_t \varphi \frac{dt}{t} \right)  = 0.
\end{equation}

Let $n > \frac{|I|}{2}$ be an integer. Similarly to \eqref{Fen3identityformula}, one has for all $x\in G$ and all $t\in (0,1)$, 
$$H_t \varphi(x) = \frac{1}{(n-1)!} \int_t^1 (v-t)^{n-1} \Delta^n H_v \varphi(x) dv + \sum_{k=0}^{n-1} \frac{1}{k!} (1-t)^k \Delta^k H_1 \varphi(x).$$
Hence, for all $x\in G$ and all $u\in (0,1)$, we have the identity
\[\begin{split}
  \int_0^u (t\Delta)^{m} H_t \varphi(x) \frac{dt}{t} 
& =  \frac{1}{(n-1)!} \int_t^1 \Delta^{n+m} H_v \varphi(x) \left(\int_{0}^{\min\{u,v\}} t^{m-1} (v-t)^{n-1} dt\right)dv \\
& \qquad  + \sum_{k=0}^{n-1} \frac{1}{k!} \Delta^{k+m} H_1 \varphi(x) \int_{0}^u t^{m-1} (1-t)^k dt.
\end{split}\]
Note that 
$$\int_{0}^{\min\{u,v\}} t^{m-1} (v-t)^{n-1} dt \lesssim u^m v^{n-1}$$
and
$$\int_{0}^u t^{m-1} (1-t)^k dt \lesssim u^m.$$
Therefore, the Schwartz seminorms of $\int_0^u (t\Delta)^{m} H_t \varphi \frac{dt}{t}$ can be estimated by
\begin{equation}\label{Fen3LPcalculus1}\begin{split}
  N_{I,c}\left( \int_0^u (t\Delta)^{m} H_t \varphi \frac{dt}{t} \right)
& \lesssim u^m \int_0^1 v^{n-1} \sup_{x\in G} e^{c|x|} |X_I \Delta^{m+n} H_v \varphi(x)| dv \\
& \qquad + u^m \sum_{k=0}^{n-1}  \sup_{x\in G} e^{c|x|} |X_I \Delta^{k+m} H_1 \varphi(x)|.
  \end{split}\end{equation}
Check then that for all $w\in (0,1]$ and all $l\in \N$, we have
\begin{equation}\label{Fen3LPcalculus2}\begin{split}
 \sup_{x\in G} e^{c|x|} |X_I \Delta^{l} H_w \varphi(x)|
& = \sup_{x\in G} e^{c|x|} \left|X_I H_w \Delta^{l}\varphi(x)\right| \\
& \leq  \sup_{x\in G} e^{c|x|}  \int_G \left|X_I h_w(y^{-1}x)\right| |\Delta^{l} \varphi(y)| dy  \\
& \lesssim  \sup_{x\in G} \int_G   e^{c|y^{-1}x|} \left|X_I h_w(y^{-1}x)\right| e^{c|y|}|\Delta^{l} \varphi(y)| dy  \\ 
& \lesssim  \left(\sup_{x\in G} \int_G   e^{c|y^{-1}x|} \left|X_I h_w(y^{-1}x)\right| dy\right)  \sum_{|I|=2l}  N_{I,c}(\varphi)
 \end{split}\end{equation}
where the third line holds because $\left\vert x\right\vert\leq \left\vert y^{-1}x\right\vert+\left\vert x\right\vert$.

\noindent However, for all $x\in G$ and all $w\in (0,1]$, Proposition \ref{Fen3htestimates} yields that, for all $x\in G$,
\begin{equation}\label{Fen3LPcalculus3}\begin{split}
   \int_G   e^{c|y^{-1}x|} \left|X_I h_w(y^{-1}x)\right| dy & \lesssim w^{-\frac{|I|}2} \frac{1}{V(\sqrt{w})} \int_G e^{c|y^{-1}x|} e^{-c'\frac{|y^{-1}x|^2}{w}} dy \\
& \lesssim w^{-\frac{|I|}2} \frac{1}{V(\sqrt{w})} \int_G  e^{-c'\frac{|y^{-1}x|^2}{2w}} dy \\
& \lesssim w^{-\frac{|I|}2}.
  \end{split}\end{equation}
By gathering the estimates \eqref{Fen3LPcalculus1}, \eqref{Fen3LPcalculus2} and \eqref{Fen3LPcalculus3}, we obtain
\[\begin{split}
   N_{I,c}\left( \int_0^u (t\Delta)^{m} H_t \varphi \frac{dt}{t} \right) 
& \lesssim u^m \left[\sum_{|I|\leq 2(m+n)}  N_{I,c}(\varphi) \right]\left[ \int_0^1 v^{n-1} v^{-\frac{|I|}{2}}dv  + \sum_{k=0}^{n-1} 1 \right] \\
& \lesssim u^m \sum_{|I|\leq 2(m+n)}  N_{I,c}(\varphi) \\
& \xrightarrow{u\to 0} 0,
  \end{split}\]
which proves \eqref{Fen3schwartzconvergence} and finishes the proof.
\end{dem}

\section{Proof of Theorem \ref{FEN3MAIN1} and of its corollaries}

\subsection{Proof of Theorem \ref{FEN3MAIN1}}

In this section, we will always assume that $\alpha\geq 0$, $p,q\in [1,+\infty]$.

\begin{prop} \label{FEN3MAIN1i1}
For all $t_1,t_0\in (0,1)$ and all integers $m > \frac{\alpha}{2}$,
 $$\|f\|_p \lesssim \|H_{t_0} f\|_p + \left( \int_0^1 \left( t^{m-\frac{\alpha}{2}} \left\| \Delta^m H_t f\right\|_p\right)^q \frac{dt}{t}\right)^{\frac{1}{q}} \qquad \forall f\in {\mathcal S}^{\prime}(G)$$
when $\alpha>0$ and
 $$\|H_{t_1}f\|_p \lesssim \|H_{t_0} f\|_p + \left( \int_0^1 \left( t^{m-\frac{\alpha}{2}} \left\| \Delta^m H_t f\right\|_p\right)^q \frac{dt}{t}\right)^{\frac{1}{q}} \qquad \forall f\in \mathcal S'(G)$$
when $\alpha \geq 0$ and $q<+\infty$, with the usual modification when $q=+\infty$.
\end{prop}

\begin{dem}
Lemma \ref{Fen3IdentityFormula} (recall that $L^p(G)\subset {\mathcal S}^{\prime}(G)$) yields the estimate
$$\|f\|_p \lesssim \int_0^{1} t^{m} \|\Delta^m H_t f\|_p \frac{dt}{t} +  \sum_{k=0}^{m-1} \|\Delta^k H_1 f\|_p.$$
However, for all $k\in \N$, $\|\Delta^k H_1 f\|_p \leq \frac{C}{(1-t_0)^k}\|H_{t_0}f\|_p$. Then, when $\alpha>0$,
\[\begin{split}
  \|f\|_p & \lesssim \int_0^{1} t^{m} \left\| \Delta^m H_t f\right\|_p \frac{dt}{t} + \|H_{t_0}f\|_p \\
& \lesssim \left(\int_0^{1}\left(t^{m-\frac{\alpha}{2}} \left\| \Delta^m H_t f\right\|_p\right)^q \frac{dt}{t}\right)^{\frac{1}{q}} \left(\int_0^{1}t^{\frac{q'\alpha}{2}} \frac{dt}{t}\right)^{\frac{1}{q'}} + \|H_{t_0}f\|_p \\
& \lesssim \left(\int_0^{1}\left(t^{m-\frac{\alpha}{2}} \left\| \Delta^m H_t f\right\|_p\right)^q \frac{dt}{t}\right)^{\frac{1}{q}} + \|H_{t_0}f\|_p,
  \end{split}\]
which prove the case $\alpha >0$. 

If $\alpha = 0$, Lemma \ref{Fen3IdentityFormula} for the integer $m+1$ implies
\[\begin{split}
  \|H_{t_1} f\|_p & \lesssim \int_0^{1} t^{m+1} \left\| \Delta^{m+1} H_{t+t_1} f\right\|_p \frac{dt}{t} + \sum_{k=0}^{m} \|\Delta^k H_{1+t_1} f\|_p \\
& \lesssim \int_0^{1} \frac{t^{m+1}}{t_1} \left\| \Delta^{m} H_{t} f\right\|_p \frac{dt}{t} + \|H_{t_0} f\|_p \\
& \lesssim \left(\int_0^{1}\left(t^{m} \left\| \Delta^m H_t f\right\|_p\right)^q \frac{dt}{t}\right)^{\frac{1}{q}} \left(\int_0^{1}\left( \frac{t}{t_1}\right)^{q'} \frac{dt}{t}\right)^{\frac{1}{q'}} + \|H_{t_0}f\|_p \\
& \lesssim \left(\int_0^{1}\left(t^{m} \left\| \Delta^m H_t f\right\|_p\right)^q \frac{dt}{t}\right)^{\frac{1}{q}} + \|H_{t_0}f\|_p.
  \end{split}\]
\end{dem}

\begin{prop} \label{FEN3MAIN1i2}
For all integers $m > \frac{\alpha}{2}$,
 $$\left( \int_0^1 \left( t^{m-\frac{\alpha}{2}} \left\| \Delta^m H_t f\right\|_p\right)^q \frac{dt}{t}\right)^{\frac{1}{q}} \lesssim \|f\|_p + \left( \int_0^1 \left( t^{m+1-\frac{\alpha}{2}} \left\| \Delta^{m+1} H_t f\right\|_p\right)^q \frac{dt}{t}\right)^{\frac{1}{q}}.$$
\end{prop}

\begin{dem}
We use Lemma \ref{Fen3IdentityFormula} and get
$$\Delta^m H_t f =  \int_0^{1} s\Delta H_{s} \Delta^m H_t f \frac{ds}{s} + H_1 \Delta^m H_t f.$$
Thus,
\[\begin{split}
&  \left( \int_0^1 \left( t^{m-\frac{\alpha}{2}} \left\| \Delta^m H_t f\right\|_p\right)^q \frac{dt}{t}\right)^{\frac{1}{q}} \\
& \qquad \leq \left( \int_0^1 \left( t^{m-\frac{\alpha}{2}} \int_0^1 \left\| \Delta^{m+1} H_{t+s} f\right\|_p ds \right)^q \frac{dt}{t}\right)^{\frac{1}{q}}  
+ \left( \int_0^1 \left( t^{m-\frac{\alpha}{2}} \left\| \Delta^m H_{t+1} f\right\|_p\right)^q \frac{dt}{t}\right)^{\frac{1}{q}} \\
& \qquad : = I_1 + I_2.
  \end{split}\]
We start with the estimate of $I_1$. One has $\Delta^{m+1}   H_{t+s} f = H_s \Delta^{m+1} H_{t} f = H_t \Delta^{m+1} H_{s} f$.
Then 
\[\begin{split}
I_1 & \leq   \left( \int_0^1 \left( t^{m-\frac{\alpha}{2}} \int_0^t \left\|\Delta^{m+1} H_{t} f \right\|_p ds \right)^q \frac{dt}{t}\right)^{\frac{1}{q}} 
+ \left( \int_0^1 \left( t^{m-\frac{\alpha}{2}} \int_t^1 \left\|\Delta^{m+1} H_{s} f \right\|_p ds \right)^q \frac{dt}{t}\right)^{\frac{1}{q}} \\
& : = II_1 + II_2.
   \end{split}\]
Notice 
$$II_1 = \left( \int_0^1 \left( t^{m+1-\frac{\alpha}{2}}\left\|\Delta^{m+1} H_{t} f \right\|_p \right)^q \frac{dt}{t}\right)^{\frac{1}{q}}$$
which is the desired estimate. As far as $II_2$ is concerned,
$$II_2 = \left( \int_0^1 \left( \int_0^1  K(s,t)g(s) \frac{ds}{s} \right)^q \frac{dt}{t} \right)^{\frac1q}$$
with $g(s) = s^{m+1-\frac{\alpha}{2}}\left\|\Delta^{m+1} H_{s} f \right\|_p$ and $K(s,t) = \left( \frac{t}{s} \right)^{m-\frac{\alpha}{2}} \1_{s\geq t}$.
Since
$$\int_0 ^1K(s,t)\frac{ds}{s} \lesssim 1 \qquad \text{ and } \qquad \int_0 ^1K(s,t)\frac{dt}{t} \lesssim 1,$$
Lemma \ref{Fen3convolutionlemma} yields then
$$II_2 \lesssim \left( \int_0^1 g(s)^q \frac{ds}{s}\right)^{\frac{1}{q}}$$
which is also the desired estimate.

It remains to estimate $I_2$. First, verify that Proposition \ref{Fen3AnalycityofHt} of $H_t$ implies $\left\| \Delta^m H_{t+1} f\right\|_p \lesssim \|f\|_p$. 
Then we obtain
$$I_2 \lesssim \|f\|_p$$
since $\int_0^1 t^{q(m-\frac{\alpha}{2})} \frac{dt}{t} < +\infty$.
\end{dem}

\begin{prop} \label{FEN3MAIN1i3}
For all integers $\beta \geq \gamma > \frac{\alpha}{2}$,
 $$\left( \int_0^1 \left( t^{\beta-\frac{\alpha}{2}} \left\| \Delta^\beta H_t f\right\|_p\right)^q \frac{dt}{t}\right)^{\frac{1}{q}} \lesssim \left( \int_0^1 \left( t^{\gamma-\frac{\alpha}{2}} \left\| \Delta^{\gamma} H_t f\right\|_p\right)^q \frac{dt}{t}\right)^{\frac{1}{q}}.$$
\end{prop}

\begin{dem}
Proposition \ref{Fen3AnalycityofHt} implies $\left\|\Delta^{\beta-\gamma} H_{\frac{t}{2}} f\right\|_p \lesssim t^{\gamma-\beta} \|f\|_p$. Then 
\[\begin{split}
   \left( \int_0^1 \left( t^{\beta-\frac{\alpha}{2}} \left\| \Delta^{\beta} H_t f\right\|_p\right)^q \frac{dt}{t}\right)^{\frac{1}{q}} 
& \lesssim \left( \int_0^1 \left( t^{\gamma-\frac{\alpha}{2}} \left\| \Delta^\gamma H_{\frac{t}{2}} f\right\|_p\right)^q \frac{dt}{t}\right)^{\frac{1}{q}} \\
& \lesssim \left( \int_0^\frac{1}{2} \left( u^{\gamma-\frac{\alpha}{2}} \left\| \Delta^\gamma H_u f\right\|_p\right)^q \frac{du}{u}\right)^{\frac{1}{q}} \\
& \quad \leq \left( \int_0^1 \left( t^{\gamma-\frac{\alpha}{2}} \left\| \Delta^\gamma H_t f\right\|_p\right)^q \frac{dt}{t}\right)^{\frac{1}{q}}.
  \end{split}\]
\end{dem}

\begin{rmq}
Propositions \ref{FEN3MAIN1i1}, \ref{FEN3MAIN1i2} and \ref{FEN3MAIN1i3} imply (i) of Theorem \ref{FEN3MAIN1}.
\end{rmq}

\begin{prop} \label{Fen3DiscreteEquivalenceA}
Let $m>\frac{\alpha}{2}$. Then
$$\|f\|_p + \left(\sum_{j\leq -1} \left[2^{j(m-\frac{\alpha}{2})} \|\Delta^m H_{2^{j}} f\|_p \right]^q \right)^{\frac1q}$$
is an equivalent norm in $B^{p,q}_\alpha(G)$.
\end{prop}

\begin{dem}
 Assertion $(i)$ in Theorem \ref{FEN3MAIN1} and the following calculus prove the equivalence of norms:
\[\begin{split}
   \left(\sum_{j\leq -1} \left[2^{j(m-\frac{\alpha}{2})} \|\Delta^m H_{2^{j}} f\|_p \right]^q \right)^{\frac1q} 
& \leq \left( \sum_{j\leq -1} \int_{2^{j-1}}^{2^j}  \left[2^{j(m-\frac{\alpha}{2})} \|\Delta^m H_{t} f\|_p  \right]^q \frac{dt}{t} \right)^\frac1q \\
& \lesssim \left( \int_{0}^{\frac12}  \left[ t^{m-\frac{\alpha}{2}} \|\Delta^m H_{t} f\|_p  \right]^q \frac{dt}{t} \right)^\frac1q \\
& \quad \leq \left( \int_{0}^{1}  \left[ t^{m-\frac{\alpha}{2}} \|\Delta^m H_{t} f\|_p  \right]^q \frac{dt}{t} \right)^\frac1q \\
& \qquad = \left( \sum_{j\leq -1} \int_{2^j}^{2^{j+1}} \left[ t^{m-\frac{\alpha}{2}} \|\Delta^m H_t f\|_p \right]^q \frac{dt}{t} \right)^\frac1q \\
& \quad \leq \left(\sum_{j\leq -1} \left[2^{j(m-\frac{\alpha}{2})} \|\Delta^m H_{2^{j}} f\|_p \right]^q \right)^{\frac1q}.
  \end{split}\]
This proves item $(ii)$ in Theorem \ref{FEN3MAIN1}. 
\end{dem}

\begin{prop} \label{Fen3DiscreteEquivalenceC}
Let $\alpha > 0$ and $l >\frac{\alpha}{2}$. Then
\begin{equation} \label{Fen3discreteL1} \|H_\frac12f\|_p + \left(\sum_{j\leq -1} \left[2^{-j\frac{\alpha}{2}} \left\|\int_{2^j}^{2^{j+1}}\left|(t\Delta)^l H_{t} f\right| \frac{dt}{t}\right\|_p \right]^q \right)^{\frac1q}\end{equation}
is an equivalent norm in $B^{p,q}_\alpha(G)$.
\end{prop}

\begin{dem}
We denote by $\|.\|_{B^{p,q}_{\alpha,1}}$ the norm defined in \eqref{Fen3discreteL1}. It is easy to check, using assertion $(i)$ in Theorem \ref{FEN3MAIN1}, the H\"older inequality and the triangle inequality, that 
$$\|f\|_{B^{p,q}_{\alpha,1}} \lesssim \|f\|_{B^{p,q}_\alpha}.$$
For the converse inequality, we proceed as follows. Fix an integer $m>\frac{\alpha}{2}$.

\begin{enumerate}
 \item {\bf Decomposition of $f$:}

The first step is to decompose $f$ as in Lemma \ref{Fen3IdentityFormula}
$$ f = \frac{1}{(l-1)!} \int_0^{1} (t\Delta)^{l} H_t f \frac{dt}{t} +  \sum_{k=0}^{l-1} \frac{1}{k!} \Delta^k H_1 f\mbox{ in }{\mathcal S}^{\prime}(G).$$

We introduce
$$f_n = -\int_{2^n}^{2^{n+1}} (t\Delta)^{l} H_t f \frac{dt}{t} dt$$
and
$$c_n = \left\| \int_{2^{n-1}}^{2^{n}} \left| (t\Delta)^{l}H_t f \right| \frac{dt}{t} \right\|_p.$$
Remark then that
$$f = \frac{1}{(l-1)!} \sum_{n=-\infty}^{-1} f_n + \sum_{k=0}^{l-1} \frac{1}{k!} \Delta^k H_1 f \qquad \text{in } \mathcal S'(G).$$

\item {\bf Estimates of $ \Delta^m H_{2^j}f_n$}

Note that
\[\begin{split} \Delta^m H_{2^j}& f_n \\
& = -\Delta^m H_{2^{n-1} + 2^j} \int_{2^{n}}^{3.2^{n-1}} (t\Delta)^l H_{t-2^{n-1}} f \frac{dt}{t}  - \Delta^m H_{2^{n} + 2^j} \int_{3.2^{n-1}}^{2^{n+1}} (t\Delta)^l H_{t-2^{n}} f \frac{dt}{t} \\
& = -\Delta^m H_{2^{n-1}+2^j} \int_{2^{n-1}}^{2^{n}} (t+2^{n-1})^l\Delta^l H_{t} f \frac{dt}{t}  - \Delta^m H_{2^{n}+2^j} \int_{2^{n-1}}^{2^{n}} (t+2^n)^l \Delta^l H_{t} f \frac{dt}{t}. \\
 \end{split}\]
Then Proposition \ref{Fen3AnalycityofHt} implies,
\begin{equation} \label{Fen3DeltaHtfn} \begin{split}
 \|\Delta^m H_{2^j}f_n \|_p & \lesssim [(2^{(n-1)}+2^j]^{-m} \left\|\int_{2^{n-1}}^{2^{n}} (t+2^{n-1})^l\Delta^l H_{t} f \frac{dt}{t} \right\|_p \\
& \qquad + [2^{n}+2^j]^{-m}\left\|\int_{2^{n-1}}^{2^{n}} (t+2^{n})^l\Delta^l H_{t} f \frac{dt}{t} \right\|_p \\
& \lesssim [2^{n}+2^j]^{-m} \left\|\int_{2^{n-1}}^{2^{n}} \left|(t\Delta)^l H_{t} f\right| \frac{dt}{t} \right\|_p.
  \end{split} \end{equation}
In other words,
\begin{equation} \label{Fen3estimdeltam}
\|\Delta^m H_{2^j} f_n\|_p \lesssim \left\{\begin{array}{ll} 2^{-nm} c_n & \text{ if } j \leq n \\ 2^{-jm} c_{n} & \text{ if } j > n \end{array} \right. .
\end{equation}

\item {\bf Estimate of $\Lambda^{p,q}_\alpha(\sum f_n)$}

As a consequence,
\[\begin{split}
   \sum_{j\leq -1}  \left[2^{j(m-\frac{\alpha}{2})} \left\|\Delta^m H_{2^{j}} \sum_{n\leq -1} f_n\right\|_p \right]^q 
& \lesssim \sum_{j\leq -1}  \left[2^{j(m-\frac{\alpha}{2})} \sum_{n\leq -1} 2^{-m\max\{j,n\}} c_n \right]^q.
\end{split} \] 
According to Lemma \ref{Fen3Bcalculus}, since $0<m-\frac\alpha 2 < m$, one has
$$\sum_{j\leq -1}  \left[2^{j(m-\frac{\alpha}{2})} \left\|\Delta^m H_{2^{j}} \sum f_n\right\|_p \right]^q  \lesssim \sum_{n=-\infty}^{-1} \left[2^{-n\frac{\alpha}{2}} c_n \right]^q.$$

\item {\bf Estimate of the remaining term}

Remark that 
$$\|f\|_{B^{p,q}_\alpha} \lesssim \|H_{t_0}f\|_p + \Lambda^{p,q}_\alpha\left(\sum f_n\right)  + \sum_{k=0}^{l-1} \Lambda^{p,q}_\alpha\left(\Delta^k H_1 f\right).$$
From the previous step and Proposition \ref{Fen3DiscreteEquivalenceA}, we proved that 
$$\Lambda^{p,q}_\alpha\left(\sum f_n\right) \lesssim \|f\|_{B^{p,q}_{\alpha,1}}.$$
In order to conclude the proof of Proposition \ref{Fen3DiscreteEquivalenceC}, it suffices then to check that 
for all $k\in \bb 0,l-1 \bn$, one has 
\begin{equation} \label{Fen3BesselNormforH1f}\|\Delta^k H_1 f\|_{B^{p,q}_\alpha} \lesssim \|f\|_{L^p}.\end{equation}
Indeed, one has for all $j \leq -1$
\[\begin{split}
   \|\Delta^m H_{2^j} \Delta^k H_1 f\|_p & = \|\Delta^{m+k} H_{1+2^j} f\|_p \\
& \lesssim \left(1+2^j\right)^{-(m+k)} \|f\|_p \\
& \lesssim \|f\|_p.
  \end{split}\]
Consequently,
\[\begin{split}
  \sum_{j\leq -1} \left[2^{j(m-\frac{\alpha}{2})} \|\Delta^m H_{2^{j}} \Delta^kH_1f\|_p \right]^q 
& \lesssim \| f\|_p^q \sum_{j\leq -1} 2^{jq(m-\frac{\alpha}{2})}  \\
& \lesssim \|f\|_p^q.
  \end{split}\]
\end{enumerate}
\end{dem}

\subsection{Proof of Theorem \ref{Fen3DiscreteEquivalenceD}}

\begin{dem} \em (Theorem \ref{Fen3DiscreteEquivalenceD}) \em

We denote by $\|.\|_{B^{p,q}_{\alpha,Xsup}}$ the norm defined in \eqref{Fen3discreteXsup}. Since 
$$ \left\Vert \Delta^{m}H_{2^j}f\right\Vert_p\leq \max_{t\in [2^j,2^{j+1}]} \sup_{\left\vert I\right\vert\leq 2m} \left\Vert X_IH_tf\right\Vert_p,$$
it is easy to check that 
$$\|f\|_{B^{p,q}_{\alpha}} \lesssim \|f\|_{B^{p,q}_{\alpha,Xsup}}.$$
For the converse inequality, it is enough to check that 
$$\|f\|_{B^{p,q}_{\alpha,Xsup}} \lesssim \|f\|_{B^{p,q}_{\alpha,1}}.$$
We proceed then as the proof of Proposition \ref{Fen3DiscreteEquivalenceC} 
since Proposition \ref{Fen3AnalycityofHt} yields 
$$
\max_{t\in [2^j,2^{j+1}]} \sup_{|I| \leq \bar m}\|X_I H_t f_n\|_p\lesssim \left\{\begin{array}{ll} 2^{-n\frac{\bar m}{2}} c_n & \text{ if } j \leq n \\ 2^{-j\frac{\bar m}{2}} c_{n} & \text{ if } j > n \end{array} \right.
$$
with a proof analogous to the one of \eqref{Fen3estimdeltam}.
\end{dem}

\subsection{Embeddings and interpolation}

\begin{dem} (of Corollary \ref{Fen3InclusionofBpqs})
The proof is analogous to the one of Proposition 2.3.2/2 in \cite{Tri1} using Proposition \ref{Fen3DiscreteEquivalenceA}. It relies on the monotonicity of $l_q$ spaces, see \cite[1.2.2/4]{Tri1}.
\end{dem}

Let us turn to interpolation properties of Besov spaces, that implies in particular Corollary \ref{Fen3InterpolationTh}.

\begin{cor} 
Let $s_0,s_1,s\geq 0$, $1\leq p_0,p_1,p,q_0,q_1,r\leq \infty$ and $\theta\in (0,1)$.

Define
$$s^* = (1-\theta)s_0 + \theta s_1,$$
$$\frac{1}{p^*} = \frac{1-\theta}{p_0} + \frac{\theta}{p_1},$$
$$\frac{1}{q^*} = \frac{1-\theta}{q_0} + \frac{\theta}{q_1}.$$
\begin{enumerate}[i.]
 \item If $s_0 \neq s_1$ then
$$(B^{p,q_0}_{s_0}, B^{p,q_1}_{s_1})_{\theta,r} = B^{p,r}_{s^*}.$$
 \item In the case where $s_0=s_1$, we have
$$(B^{p,q_0}_{s}, B^{p,q_1}_{s})_{\theta,q^*} = B^{p,q^*}_{s}.$$
 \item If $p^* = q^* : = r$,
$$(B^{p_0,q_0}_{s_0}, B^{p_1,q_1}_{s_1})_{\theta,r} = B^{r,r}_{s^*}.$$
 \item If $s_0 \neq s_1$,
$$(B^{p_0,q_0}_{s_0}, B^{p_1,q_1}_{s_1})_{[\theta]} = B^{p^*,q^*}_{s^*}.$$
\end{enumerate}
\end{cor}

\begin{dem}
The proof is inspired by \cite[Theorem 6.4.3]{BL}. 

Recall (see Definition 6.4.1 in \cite{BL}) that a space $B$ is called a retract of $A$ if there exists two bounded linear operators $\mathcal J:B\to A$ and $\mathcal P:A\to B$ such that $\mathcal P \circ \mathcal J$ is the identity on $B$.

\noindent Therefore, we just need to prove that the spaces $B^{p,q}_\alpha$ are retracts of $l^\alpha_q(L^p)$ where, for any Banach space $A$ (see paragraph 5.6 in \cite{BL}),
$$l^\alpha_q(A) = \left\{ u \in A^{\Z_-}, \, \|u\|_{\ell^\alpha_q(A)} := \left(\sum_{j\leq 0} \left[2^{-j\frac{\alpha}{2}} \|u_j\|_{A}\right]^q \right)^\frac1q < +\infty \right\}.$$
Then interpolation on the spaces $l^\alpha_q(L^p)$ (see \cite{BL}, Theorems 5.6.1, 5.6.2 and 5.6.3) provides the result.
Note the weight appearing $l^\alpha_q(A)$ is $2^{-j\frac{\alpha}{2}}$ (and not $2^{j\frac{\alpha}{2}}$) because we sum on \em negative \em integers.

\ms

Fix $m>\frac{\alpha}2$. Define the functional $\mathcal J$  by $\mathcal J f = \left((\mathcal J f)_j\right)_{j\leq 0}$ where
$$(\mathcal J f)_j = 2^{jm} \Delta^m H_{2^{j-1}} f$$
if $j\leq -1$ and
$$(\mathcal J f)_0 = H_\frac12 f.$$
Moreover, define $\mathcal P$ on $l^\alpha_q(L^p)$ by
$$\mathcal P u = \sum_{k=0}^{2m-1} \frac{1}{k!} \Delta^k H_\frac12 u_0 + \frac{1}{(2m-1)!}\sum_{j\leq -1} 2^{-jm} \int_{2^j}^{2^{j+1}} t^{2m}\Delta^m H_{t-2^{j-1}} u_j \frac{dt}{t}.$$
We will see below that $\mathcal P$ is well-defined on $l^\alpha_q(L^p)$. Proposition \ref{Fen3DiscreteEquivalenceA} implies immediately that $\mathcal J$ is bounded from $B^{p,q}_\alpha$ to $l^\alpha_q(L^p)$.
Moreover, Lemma \ref{Fen3IdentityFormula} easily provides that
$$\mathcal P \circ \mathcal J  = Id_{B^{p,q}_\alpha}.$$
It remains to verify that $\mathcal P$ is a bounded linear operator from $l^\alpha_q(L^p)$ to $B^{p,q}_\alpha$. The proof is similar to the one of Proposition \ref{Fen3DiscreteEquivalenceC}.
Indeed, proceeding as the fourth step of Proposition \ref{Fen3DiscreteEquivalenceC}, one gets 
$$\left\|\sum_{k=0}^{2m-1} \frac{1}{k!} \Delta^k H_\frac{1}{2} u_0\right\|_{B^{p,q}_\alpha} \lesssim \|u_0\|_p.$$
It is plain to see that
\[\begin{split}
   \left\|H_\frac12 \sum_{j\leq -1} 2^{-jm} \int_{2^j}^{2^{j+1}} t^{2m}\Delta^m H_{t-2^{j-1}} u_j \frac{dt}{t} \right\|_p 
& \leq \sum_{j\leq -1} 2^{-jm} \int_{2^j}^{2^{j+1}} t^{2m} \left\|H_\frac12 \Delta^m H_{t-2^{j-1}} u_j  \right\|_p \frac{dt}{t} \\
& \leq \sum_{j\leq -1} 2^{-jm} \int_{2^j}^{2^{j+1}} t^{2m} \left\|H_\frac12 \Delta^m u_j  \right\|_p \frac{dt}{t} \\
& \lesssim \sum_{j\leq -1} 2^{jm} \|u_j\|_p \\
& \lesssim \|u\|_{l^\alpha_q(L^p)}.
  \end{split}\]

Then the proof of the boundedness of $\mathcal P$ is reduced to the one of
\begin{equation} \label{Fen3lastassumption}
I: = \left( \sum_{k\leq -1} \left( 2^{k(m-\frac{\alpha}{2})} \left\|\Delta^m H_{2^k} \sum_{j\leq -1} 2^{-jm} \int_{2^j}^{2^{j+1}} t^{2m}\Delta^m H_{t-2^{j-1}} u_j \frac{dt}{t} \right\|_p \right)^q\right)^\frac1q \lesssim \|u\|_{l^\alpha_q(L^p)}.
\end{equation}
Indeed,
\[\begin{split}
   I^q & \lesssim \sum_{k\leq -1} \left( 2^{k(m-\frac{\alpha}{2})} \left\| \sum_{j\leq -1} 2^{-jm} \int_{2^j}^{2^{j+1}} (t\Delta)^{2m} H_{t-2^{j-1}+2^k} u_j \frac{dt}{t} \right\|_p \right)^q \\
& \lesssim \sum_{k\leq -1} \left( 2^{k(m-\frac{\alpha}{2})}  \sum_{j\leq -1} 2^{jm} \left\|\Delta^{2m} H_{2^{j-1}+2^k} u_j \right\|_p \right)^q \\
& \lesssim \sum_{k\leq -1} \left( 2^{k(m-\frac{\alpha}{2})}  \sum_{j\leq -1} \frac{2^{jm}}{(2^j+2^k)^{2m}} \left\| u_j \right\|_p \right)^q \\
& \lesssim \sum_{k\leq -1} \left( 2^{k(m-\frac{\alpha}{2})}  \sum_{j\leq -1} 2^{-2m\max\{j,k\}} 2^{jm}\left\| u_j \right\|_p \right)^q
  \end{split}\]

Check that $0<m-\frac{\alpha}{2}<2m$. Thus, Lemma \ref{Fen3Bcalculus} yields
$$I^q \lesssim \sum_{j\leq -1} \left[ 2^{-j(\frac{\alpha}{2}+m)} 2^{jm}\|u_j\|_p\right]^q \leq \|u\|_{l^\alpha_q(L^p)}^q,$$
which proves \eqref{Fen3lastassumption} and thus concludes the proof.
\end{dem}

\section{Algebra under pointwise product - \\ Theorem \ref{FEN3MAIN2}}

We want to introduce some paraproducts. The idea of paraproducts goes back to \cite{Bony}. 
The term ``paraproducts'' is used to a denotes some non-commutative bilinear forms $\Lambda_i$ such that $fg = \sum \Lambda_i(f,g)$. 
They are introduced in some cases, where the bilinear forms $\Lambda_i$  are easier to handle than the pointwise product.

\noindent In the context of doubling spaces, a definition of paraproducts is given in \cite{Bernicot,FreyThesis}.
We need to slightly modify the definition in \cite{Bernicot} to adapt them to non-doubling spaces.

\ms

For all $t>0$, define 
$$\phi_t(\Delta) = -\sum_{k=0}^{m-1} \frac{1}{k!} (t\Delta)^k H_t,$$
and observe that the derivative of $t\mapsto \phi_t(\Delta)$ is given by
$$\phi'_t(\Delta) = \frac{1}{(m-1)!} \frac 1t (t\Delta)^m H_t:=\frac 1t \psi_t(\Delta).
$$
\begin{rmq}
Even if $\phi_t$ actually depends on $m$, we do not indicate this dependence explicitly.
\end{rmq}
Recall that Lemma \ref{Fen3IdentityFormula} provides the identity
\begin{equation} \label{Fen3IdentityFormula2} f = \int_0^1 \psi_t(\Delta) f \frac{dt}{t} - \phi_1(\Delta) f \qquad \text{in } \mathcal S'(G).\end{equation}

\begin{prop} \label{Fen3ParaproductFormula}
Let $p,q,r\in [1,+\infty]$ such that $\frac{1}{r}: = \frac{1}{p} + \frac{1}{q} \leq 1$. Let $(f,g) \in L^p(G)\times L^q(G)$.

One has the formula
$$fg = \Pi_f(g) + \Pi_g(f) + \Pi(f,g) - \phi_1(\Delta)[\phi_1(\Delta) f\cdot \phi_1(\Delta) g] \qquad \text{in } \mathcal S'(G),$$
where
$$\Pi_f(g) = \int_0^1 \phi_t(\Delta)[\psi_t(\Delta) f \cdot \phi_t(\Delta)g] \frac{dt}{t}$$
and
$$\Pi(f,g) = \int_0^1 \psi_t(\Delta)[\phi_t(\Delta) f \cdot \phi_t(\Delta)g] \frac{dt}{t}.$$
\end{prop}

\begin{dem}
Since $fg \in L^r \subset \mathcal S'(G)$, the formula \eqref{Fen3IdentityFormula2} provides in $\mathcal S'(G)$
\begin{equation}
 [f\cdot g] = \int_0^1 \psi_t(\Delta) [f\cdot g] \frac{dt}{t} - \phi_1(\Delta) [f\cdot g]. 
\end{equation}
We can use again twice (one for $f$ and one for $g$) the identity \eqref{Fen3IdentityFormula2} to get
\begin{equation}
\begin{split}
 [f\cdot g] & = \int_0^1 \psi_t(\Delta) \left[\left\{\int_0^1 \psi_u(\Delta) f \frac{du}{u} - \phi_1(\Delta) f\right\}\cdot \left\{\int_0^1 \psi_v(\Delta) g \frac{dv}{v} - \phi_1(\Delta) g\right\}\right] \frac{dt}{t} \\
& \qquad - \phi_1(\Delta) \left[\left\{\int_0^1 \psi_u(\Delta) f \frac{du}{u} - \phi_1(\Delta) f\right\}\cdot \left\{\int_0^1 \psi_v(\Delta) g \frac{dv}{v} - \phi_1(\Delta) g\right\}\right] \\
& = \int_0^1\int_0^1\int_0^1 \psi_t(\Delta) [\psi_u(\Delta) f \cdot \psi_v(\Delta) g] \frac{dt\, du\, dv}{tuv}  \\
& \qquad - \int_0^1\int_0^1 \psi_t(\Delta) [\phi_1(\Delta) f \cdot \psi_v(\Delta) g] \frac{dt\, dv}{tv} - \int_0^1\int_0^1 \psi_t(\Delta) [\psi_u(\Delta) f \cdot \phi_1(\Delta) g] \frac{dt\, du}{tu} \\
& \qquad \qquad \qquad \qquad \qquad \qquad \qquad \qquad \qquad \qquad \qquad \qquad -  \int_0^1\int_0^1 \phi_1(\Delta) [\psi_u(\Delta) f \cdot \psi_v(\Delta) g] \frac{du\, dv}{uv} \\
& \qquad + \int_0^1 \psi_t(\Delta) [\phi_1(\Delta) f \cdot \phi_1(\Delta) g] \frac{dt}{t} +  \int_0^1 \phi_1(\Delta) [\psi_u(\Delta) f \cdot \phi_1(\Delta) g] \frac{du}{u} \\
& \qquad \qquad \qquad \qquad \qquad \qquad \qquad \qquad \qquad \qquad \qquad \qquad + \int_0^1 \phi_1(\Delta) [\phi_1(\Delta) f \cdot \psi_v(\Delta) g] \frac{dv}{v} \\
& \qquad - \phi_1(\Delta) [\phi_1(\Delta)\cdot \phi_1(\Delta)] \\
& \qquad := R(f,g) + \int_0^1\int_0^1\int_0^1 \psi_t(\Delta) [\psi_u(\Delta) f \cdot \psi_v(\Delta) g] \frac{dt\, du\, dv}{tuv} - \phi_1(\Delta) [\phi_1(\Delta)\cdot \phi_1(\Delta)].
\end{split}
\end{equation}
The domain $[0,1]^3$ can be divided in the subsets $D(t,u,v)$, $D(u,t,v)$ and $D(v,u,t)$ where $D(a,b,c) = \{(a,b,c)\in [0,1]^3, \, a<\min\{b,c\}\}$. 
Consequently,
\begin{equation} \label{compPar}
\begin{split}
\int_0^1\int_0^1\int_0^1 &  \psi_t(\Delta) [\psi_u(\Delta) f \cdot \psi_v(\Delta) g] \frac{dt\, du\, dv}{tuv} \\
&  = \int_0^1\int_t^1\int_t^1 \psi_t(\Delta) [\psi_u(\Delta) f \cdot \psi_v(\Delta) g] \frac{dt\, du\, dv}{tuv} 
+ \int_0^1\int_u^1\int_u^1 \psi_t(\Delta) [\psi_u(\Delta) f \cdot \psi_v(\Delta) g] \frac{du\, dt\, dv}{utv} \\
& \qquad \qquad \qquad \qquad \qquad \qquad \qquad \qquad + \int_0^1\int_v^1\int_v^1 \psi_t(\Delta) [\psi_u(\Delta) f \cdot \psi_v(\Delta) g] \frac{dv\, du\, dt}{vut} \\
& = \int_0^1 \psi_t(\Delta) [\{\phi_1(\Delta) f-\phi_t(\Delta) f\} \cdot \{\phi_1(\Delta) g-\phi_t(\Delta) g\}] \frac{dt}{t} \\
&  \qquad \qquad \qquad \qquad \qquad \qquad \qquad + \int_0^1 \{\phi_1(\Delta) - \phi_u(\Delta) \} [\psi_u(\Delta) f \cdot \{\phi_1(\Delta) g-\phi_u(\Delta) g\}] \frac{du}{u} \\
&  \qquad \qquad \qquad \qquad \qquad \qquad \qquad + \int_0^1 \{\phi_1(\Delta) - \phi_v(\Delta) \} [\{\phi_1(\Delta) f-\phi_v(\Delta) f\} \cdot \psi_v(\Delta) g] \frac{dv}{v} \\
& := S(f,g) + \Pi_f(g) + \Pi_g(f) + \Pi(f,g).
\end{split}
\end{equation}
It remains to check that $R(f,g) + S(f,g) = 0$. 
This identity, that can be proven with similar computations as \eqref{compPar}, is left to the reader.
\end{dem}

\begin{prop} \label{Fen3Paraproduct1}
Let $G$ be a unimodular Lie group. Let $\alpha > 0$ and $p,p_1,p_2,q\in [1,+\infty]$ such that
$$\frac{1}{p_1} + \frac{1}{p_2} = \frac{1}{p}.$$ 
Then for all $f\in B^{p_1,q}_\alpha$ and all $g \in L^{p_2}$, one has
$$\Lambda^{p,q}_\alpha[\Pi_f(g)] \lesssim \|f\|_{B^{p_1,q}_\alpha}\|g\|_{L^{p_2}}.$$ 
\end{prop}

\begin{dem}
 Let $m > \frac{\alpha}{2}$ and $j\leq -1$. Notice that, for all $u\in (0,1)$,
$$ \left\| \Delta^ m H_{u} \Pi_f(g) \right\|_p \leq  \int_0^1 \left\|\Delta^m H_u \phi_t(\Delta)[\psi_t(\Delta) f \cdot \phi_t(\Delta)g] \right\|_p \frac{dt}{t}.$$
Remark that $$\|\phi_t(\Delta) h\|_r \lesssim \|H_{\frac{t}{2}} h\|_r$$ for all $r\in [1,+\infty]$ and all $h \in L^r$. As a consequence,
\[\begin{split}
   \left\|\Delta^m H_{u} \phi_t(\Delta)[\psi_t(\Delta) f \cdot \phi_t(\Delta)g] \right\|_p 
& = \left\|\phi_t(\Delta)\Delta^m H_{u} [\psi_t(\Delta) f \cdot \phi_t(\Delta)g] \right\|_p \\
& \lesssim   \left\|\Delta^m H_{\frac t2+u} [\psi_t(\Delta) f \cdot \phi_t(\Delta)g] \right\|_p\\
& \lesssim \left(\frac t2+u\right)^{-m} \left\|\psi_t(\Delta) f \cdot \phi_t(\Delta)g \right\|_p\\
& \lesssim \min\left\{t^{-m}, u^{-m}\right\}  \left\|\psi_t(\Delta) f \cdot \phi_t(\Delta)g \right\|_p \\
& \lesssim \min\left\{t^{-m}, u^{-m}\right\}  \left\|\psi_t(\Delta) f\right\|_{p_1} \left\| \phi_t(\Delta)g \right\|_{p_2} \\
& \lesssim \min\left\{t^{-m}, u^{-m}\right\}  \left\|(t\Delta)^m H_t f\right\|_{p_1} \left\|g \right\|_{p_2}.
  \end{split}\]
We deduce then
\[\begin{split}
   \Lambda^{p,q}_\alpha & [\Pi_f (g)]^q \\
& \lesssim \left\|g \right\|_{p_2} ^q\int_0^1 \left( u^{m-\frac{\alpha}{2}} \int_{0}^{1} \left(\max\{u,t\}\right)^{-m} \left\|(t\Delta)^m H_t f\right\|_{p_1} \frac{dt}{t}  \right)^q \frac{du}{u} \\
& \lesssim \left\|g \right\|_{p_2} ^q \sum_{j\leq -1} \left( 2^{j(m-\frac{\alpha}{2})} \sum_{n=-\infty}^{-1} 2^{-m\max\{j,n\}} \left\|(2^n\Delta)^m H_{2^n} f\right\|_{p_1}   \right)^q \\
& \lesssim \|g\|_{L^{p_2}} \left(\sum_{n\leq -1} 2^{-n\frac{\alpha}{2}q} 2^{nmq}\|\Delta^m H_{2^n} f\|_{p_1}^q \right)^\frac1q,
\end{split}\]
where we used Lemma \ref{Fen3Bcalculus} for the last line.
As a consequence, we obtain if $\alpha \in (0,2m)$,
\[\begin{split}
 \Lambda^{p,q}_\alpha[\Pi_f(g)] & \lesssim \|g\|_{L^{p_2}} \left(\sum_{n\leq -1} 2^{nq(m-\frac{\alpha}{2})} \|\Delta^m H_{2^n} f\|_{p_1}^q \right)^\frac1q \\
& \lesssim \|g\|_{L^{p_2}} \|f\|_{B^{p_1,q}_\alpha}  
  \end{split}\]
where we used Proposition \ref{Fen3DiscreteEquivalenceA} for the last line.
\end{dem}

\begin{prop} \label{Fen3Paraproduct2}
Let $G$ be a unimodular Lie group. Let $\alpha > 0$ and $p,p_1,p_2,p_3,p_4,q\in [1,+\infty]$ such that
$$\frac{1}{p_1} + \frac{1}{p_2} = \frac{1}{p_3} + \frac{1}{p_4} = \frac{1}{p}.$$ 
Then for all $f\in B^{p_1,q}_\alpha \cap L^{p_3}$ and all $g \in B^{p_4,q}_\alpha \cap L^{p_2}$, one has
$$\Lambda^{p,q}_\alpha[\Pi(f,g)] \lesssim \|f\|_{B^{p_1,q}_\alpha}\|g\|_{L^{p_2}} + \|f\|_{L^{p_3}} \|g\|_{B^{p_4,q}_\alpha}.$$
\end{prop}

\begin{dem}
Notice first that
$$ \left\| \Delta^ m H_{u} \Pi(f,g) \right\|_p \leq  \int_0^1  \left\|\Delta^m H_{u} H_t (t\Delta)^m[\phi_t(\Delta) f \cdot \phi_t(\Delta)g] \right\|_p \frac{dt}{t}.$$
Let us recall then that $X_i(f \cdot g) = f \cdot X_ig + X_i f   \cdot g$. 
Consequently, since $\Delta = \sum_{i=1}^k X_i^2$, one has
$$\|\Delta^m [f\cdot g]\|_p \lesssim \|\Delta^m f  \cdot  g\|_p + \| f  \cdot \Delta^m g\|_p+\sum_{k=1}^{2m-1} \sup_{|I_1| = k} \sup_{|I_2| = 2m-k} \|X_{I_1} f  \cdot X_{I_2} g\|_p.$$
In the following computations, $(Y_{I_1},Z_{I_2})$ denotes the couple $(X_{I_1},X_{I_2})$ if $|I_1|\neq 0$ and  $|I_2|\neq 0$, $(\Delta^{|I_1|/2},I)$ if $|I_2| = 0$ and $(I,\Delta^{|I_2|/2})$ if $|I_1| = 0$.
With these notations, one has
\[\begin{split}
   \|\Delta^m H_{u+t} (t\Delta)^m & [\phi_t(\Delta) f \cdot \phi_t(\Delta)g] \|_p \\
& \lesssim \min\left\{t^{-m}, u^{-m}\right\}  \left\|(t\Delta)^m[\phi_t(\Delta) f \cdot \phi_t(\Delta)g]\right\|_p \\
& \lesssim \min\left\{t^{-m}, u^{-m}\right\}  \sum_{k,l=0}^{m-1} \left\|(t\Delta)^m[(t\Delta)^k H_t f \cdot (t\Delta)^l H_t g]\right\|_p \\
& \lesssim \min\left\{t^{-m}, u^{-m}\right\}  \sum_{k,l=0}^{m-1} \sum_{i=0}^{2m} t^m \sup_{|I_1| = i} \sup_{|I_2|=2m-i} \left\| Y_{I_1}(t\Delta)^kH_t f \cdot Z_{I_2}(t\Delta)^l H_tg \right\|_p \\
& = \min\left\{t^{-m}, u^{-m}\right\}  \sum_{k,l=0}^{m-1} \sum_{i=0}^{2m} t^{m+k+l} \sup_{|I_1| = i} \sup_{|I_2|=2m-i} \left\| Y_{I_1}\Delta^kH_t f \cdot Z_{I_2}\Delta^l H_tg \right\|_p \\
& \lesssim \min\left\{t^{-m}, u^{-m}\right\}  \sum_{k,l=0}^{m-1} \sum_{i=0}^{2m} t^{m+k+l}\sup_{|I_1|= i+2k} \sup_{|I_2|= 2m+2l-i} \left\| Y_{I_1}H_t f \cdot Z_{I_2} H_tg \right\|_p \\
& \lesssim \min\left\{t^{-m}, u^{-m}\right\}  \sum_{\begin{subarray}{c} 2m\leq k+l \leq 6m-4 \\ k+l \text{ even} \end{subarray}} t^{\frac{k+l}{2}} \sup_{|I_1|= k} \sup_{|I_2|= l}   \left\| Y_{I_1}H_t f \cdot Z_{I_2} H_tg \right\|_p.   
\end{split}\]

Setting $\ds c_n = \sum_{\begin{subarray}{c} 2m\leq k+l \leq 6m-4 \\ k+l \text{ even} \end{subarray}} \int_{2^n}^{2^{n+1}} t^{\frac{k+l}{2}} \sup_{|I_1|= k} \sup_{|I_2|= l}   \left\| Y_{I_1}H_t f \cdot Z_{I_2} H_tg \right\|_p \frac{dt}{t}$, one has
\[\begin{split}
   \Lambda^{p,q}_\alpha & [\Pi(f,g)]^q \\
& \lesssim \int_0^1 \left( u^{m-\frac{\alpha}{2}} \int_{0}^1 \|\Delta^m H_{u+t} (t\Delta)^m  [\phi_t(\Delta) f \cdot \phi_t(\Delta)g] \|_p \frac{dt}{t} \right)^q \frac{du}{u} \\
& \lesssim \int_0^1 \left( u^{m-\frac{\alpha}{2}} \int_{0}^1 \min\left\{t^{-m}, u^{-m}\right\} \sum_{\begin{subarray}{c} 2m\leq k+l \leq 6m-4 \\ k+l \text{ even} \end{subarray}} t^{\frac{k+l}{2}} \sup_{|I_1|= k} \sup_{|I_2|= l}   \left\| Y_{I_1}H_t f \cdot Z_{I_2} H_tg \right\|_p \frac{dt}{t} \right)^q \frac{du}{u} \\
& \lesssim \sum_{j=-\infty}^{-1} \left[2^{j(m-\frac{\alpha}{2})} \sum_{n=-\infty}^{-1} 2^{-m\max\{n,j\}} c_n\right]^q \\
& \lesssim \sum_{n\leq -1} 2^{-nq\frac{\alpha}{2}} c_n^q
\end{split}\]
where the last line is a consequence of Lemma \ref{Fen3Bcalculus}, since $0<m-\frac\alpha 2<m$.

It remains to prove that for any couple $(k,l) \in \N^2$ satisfying $6m-4\geq k+l \geq 2m$ and $k+l$ even, we have 
\begin{equation} \begin{split}
 T & : = \left(\sum_{n\leq -1} 2^{-nq\frac{\alpha}{2}} \left( \int_{2^n}^{2^{n+1}}  t^{\frac{k+l}{2}} \sup_{|I_1|= k} \sup_{|I_2|= l}   \left\| Y_{I_1}H_t f \cdot Z_{I_2} H_tg \right\|_p \frac{dt}{t} \right)^q \right)^\frac1q \\
& \lesssim \|f\|_{B^{p_1,q}_\alpha}\|g\|_{L^{p_2}} + \|f\|_{L^{p_3}} \|g\|_{B^{p_4,q}_\alpha}.
\end{split}\end{equation}

\begin{enumerate}
 \item {\bf If $k=0$ or $l=0$:}

Since $k$ and $l$ play symmetric roles, we can assume without loss of generality that $l=0$. In this case, $k$ is even and if $k=2k'$, 
\[\begin{split} 
  \sup_{|I_1|= k} \sup_{|I_2|= 0}   \left\| Y_{I_1}H_t f \cdot Z_{I_2} H_tg \right\|_p &  =  \left\|\Delta^{k'}H_t f \cdot H_tg \right\|_p \\
& \leq  \left\| \Delta^{k'}H_t f\|_{p_1}  \|H_tg \right\|_{p_2} \\
& \leq  \left\| \Delta^{k'} H_t f\right\|_{p_1}  \left\|g \right\|_{p_2}.
 \end{split}\]
Therefore,
\[\begin{split} 
T &  \leq \|g\|_{L^{p_2}} \left(\sum_{n\leq -1} 2^{-nq\frac{\alpha}{2}} \left( \int_{2^n}^{2^{n+1}}  t^{k'}  \left\| \Delta^{k'}H_t f \right\|_{p_1} \frac{dt}{t} \right)^q \right)^\frac1q \\
& \lesssim \|g\|_{L^{p_2}} \|f\|_{B^{p_1,q}_\alpha}
\end{split} \]
where the second line is due to the fact that $k' \geq m > \frac{\alpha}{2}$.

 \item {\bf If $k\geq 1$ and $l\geq 1$:}

Define $\alpha_1,\alpha_2,r_1,r_2,q_1$ and $q_2$ by
$$\alpha_1 = \frac{k}{k+l}\alpha \qquad , \qquad \alpha_2 =\frac{l}{k+l}\alpha,$$
$$\frac{k+l}{r_1} = \frac{k}{p_1} + \frac{l}{p_3} \qquad, \qquad \frac{k+l\label{Fen3Theo3}}{r_2} = \frac{k}{p_2} + \frac{l}{p_4},$$
$$\frac{k+l}{q_1} = \frac{k}{q} \qquad , \qquad \frac{k+l}{q_2} = \frac{l}{q}.$$
In this case, notice that $k > \alpha_1$ and $l > \alpha_2$. One has then
\[\begin{split}
  \sup_{|I_1|= k} \sup_{|I_2|= l}   \left\| X_{I_1}H_t f \cdot X_{I_2} H_tg \right\|_p 
& \leq  \sup_{|I_1|= k} \left\| X_{I_1}H_t f \right\|_{r_1}  \sup_{|I_2|= l}\left\|X_{I_2} H_tg \right\|_{r_2} 
  \end{split}\]
and thus Hölder inequality provides
\[\begin{split}
   T & \leq \left(\sum_{n\leq -1} \left( 2^{n\frac{k-\alpha_1}{2}} \max_{t\in [2^j,2^{j+1}]}   \sup_{|I_1|= k}  \left\| X_{I_1}H_t f \right\|_{r_1}  \right)^{q_1} \right)^\frac1{q_1}  \\
& \qquad \left(\sum_{n\leq -1} \left( 2^{n\frac{l-\alpha_2}{2}} \max_{t\in [2^j,2^{j+1}]}  \sup_{|I_2|= l}  \left\| X_{I_2}H_t g \right\|_{r_2}  \right)^{q_2} \right)^\frac1{q_2} \\
& \lesssim \|f\|_{B^{r_1,q_1}_{\alpha_1}} \|g\|_{B^{r_2,q_2}_{\alpha_2}}
  \end{split}\]
where the second line is due to Theorem \ref{Fen3DiscreteEquivalenceD}.

Let $\theta = \frac{k}{k+l}$. Complex interpolation (Corollary \ref{Fen3InterpolationTh}) provides
$$(B^{p_3,\infty}_0, B^{p_1,q}_{\alpha})_{[\theta]} = B^{r_1,q_1}_{\alpha_1}$$
and
$$(B^{p_4,q}_{\alpha}, B^{p_2,\infty}_0)_{[\theta]} = B^{r_2,q_2}_{\alpha_2}.$$
Remark also that $L^s(G)$ is continuously embedded in $B^{s,\infty}_0(G)$ (this can be easily seen from the definition of Besov spaces). As a consequence,
\[\begin{split}
   T & \lesssim \|f\|_{B^{r_1,q_1}_{\alpha_1}} \|g\|_{B^{r_2,q_2}_{\alpha_2}} \\
& \lesssim \|f\|_{L^{p_3}}^\theta \|f\|_{B^{p_1,q}_\alpha}^{1-\theta} \|g\|_{B^{p_4,q}_\alpha}^{\theta}  \|g\|_{L^{p_2}}^{1-\theta} \\
& \lesssim \|f\|_{L^{p_3}} \|g\|_{B^{p_4,q}_\alpha} + \|f\|_{B^{p_1,q}_\alpha}  \|g\|_{L^{p_2}}
  \end{split}\]
which is the desired conclusion.
\end{enumerate}
\end{dem}

Let us now prove Theorem \ref{FEN3MAIN2}

\begin{dem}
With the use of Propositions \ref{Fen3ParaproductFormula}, \ref{Fen3Paraproduct1} and \ref{Fen3Paraproduct2}, it remains to check that
\begin{equation} \label{Fen3restant1} \|H_\frac{1}{2}[f\cdot g]\|_{L^p} \lesssim \|f\|_{B^{p_1,q}_\alpha}\|g\|_{L^{p_2}} + \|f\|_{L^{p_3}} \|g\|_{B^{p_4,q}_\alpha}\end{equation}
and
\begin{equation} \label{Fen3restant2} \|\phi_1(\Delta)[\phi_1(\Delta)f \cdot \phi_1(\Delta) g\|_{B^{p,q}_\alpha} \lesssim \|f\|_{B^{p_1,q}_\alpha}\|g\|_{L^{p_2}} + \|f\|_{L^{p_3}} \|g\|_{B^{p_4,q}_\alpha}.\end{equation}
The inequality \eqref{Fen3restant1} is easy to check. By Proposition \ref{FEN3MAIN1i1}, one has
$$\|H_\frac{1}{2}[f\cdot g]\|_{L^p} \leq \|f\cdot g\|_p \leq \|f\|_{p_1}\|g\|_{p_2} \leq \|f\|_{B^{p_1,q}_\alpha} \|g\|_{L^{p_2}}.$$
For \eqref{Fen3restant2}, recall that \eqref{Fen3BesselNormforH1f} implies
\[\begin{split}
   \|\phi_1(\Delta)[\phi_1(\Delta)f \cdot \phi_1(\Delta) g]\|_{B^{p,q}_\alpha} & \lesssim \|\phi_1(\Delta)f \cdot \phi_1(\Delta) g\|_{L^p} \\
& \lesssim \|\phi_1(\Delta) f\|_{p_1} \|\phi_1(\Delta) g\|_{p_2} \\
& \lesssim \|f\|_{B^{p_1,q}_\alpha} \|g\|_{L^{p_2}}.
  \end{split}\]
\end{dem}

\section{Other characterizations of Besov spaces}

\subsection{Characterization by differences of functions - Theorem \ref{Fen3Main3}}

\begin{lem} \label{Fen3LambdaBoundedByL}
Let $p,q\in [1,+\infty]$ and $\alpha>0$. There exists $c>0$ such that, for all $f\in L^p(G)$,
$$\Lambda^{p,q}_\alpha(f) \lesssim \left( \int_G \left( \frac{\|\nabla_y f\|_p e^{-c|y|^2}}{|y|^\alpha} \right)^q \frac{dy}{V(|y|)} \right)^{\frac{1}{q}}.$$ 
\end{lem}

\begin{dem} 

Since $ \int_G \frac{\dr h_t}{\dr t}(y) dx = 0$,
\begin{equation}\label{Fen3Nablayappears}\begin{split}
   \frac{\dr H_t}{\dr t} f(x) & = \int_G \frac{\dr h_t}{\dr t}(y) f(xy) dy \\
& = \int_G  \frac{\dr h_t}{\dr t}(y) [f(xy) - f(x)]dy \\
& = \int_G  \frac{\dr h_t}{\dr t}(y) \nabla_y f(x)dy.
  \end{split}\end{equation}
Consequently,
$$\left\| \frac{\dr H_t}{\dr t} f\right\|_p \leq \int_G \left|\frac{\dr h_t}{\dr t}(y)\right| \|\nabla_y f\|_p dy.$$
Proposition \ref{Fen3htestimates} provides
\[\begin{split}
   \Lambda^{p,q}_\alpha(f) & \lesssim \left( \int_0^1 \left( t^{1-\frac{\alpha}{2}} \int_G \frac{1}{tV(\sqrt{t})} e^{-c\frac{|y|^2}{t}}  \|\nabla_y f\|_p dy \right)^q \frac{dt}{t}\right)^{\frac{1}{q}} \\
& \lesssim \left( \int_0^1 \left( t^{1-\frac{\alpha}{2}} \int_G \frac{1}{tV(\sqrt{t})} e^{-c'\frac{|y|^2}{t}}  \|\nabla_y f\|_p e^{-c'\frac{|y|^2}t}dy \right)^q \frac{dt}{t}\right)^{\frac{1}{q}}\\   
& \lesssim \left( \int_0^1 \left( t^{1-\frac{\alpha}{2}} \int_G \frac{1}{tV(\sqrt{t})} e^{-c'\frac{|y|^2}{t}}  \|\nabla_y f\|_p e^{-c'|y|^2}dy \right)^q \frac{dt}{t}\right)^{\frac{1}{q}} \\
& \quad = \left( \int_0^1 \left( \int_G K(t,y)  g(y) \frac{dy}{V(|y|)} \right)^q \frac{dt}{t}\right)^{\frac{1}{q}}
  \end{split}\]
with $c' = \frac{c}{2}$, $g(y) = \frac{\|\nabla_y f\|_p}{|y|^\alpha} e^{-c'|y|^2}$ and $K(t,y) = \frac{V(|y|)}{V(\sqrt{t})} \left(\frac{|y|^2}{t} \right)^\frac{\alpha}{2} e^{-c'\frac{|y|^2}{t}}$  (note that we used the fact that $t\in (0,1)$ in the third line).  
Lemma \ref{Fen3convolutionlemma}  and Proposition \ref{Fen3convolutionbyhtbis} imply  then
\[\begin{split}
   \Lambda^{p,q}_\alpha(f) & \lesssim \left( \int_G |g(y)|^q \frac{dy}{V(|y|)} \right)^{\frac{1}{q}} \\
& \quad = \left( \int_G \left( \frac{\|\nabla_y f\|_p e^{-c|y|^2}}{|y|^\alpha} \right)^q \frac{dy}{V(|y|)} \right)^{\frac{1}{q}}.
  \end{split}\]
\end{dem}

\begin{prop}
Let $p,q\in [1,+\infty]$ and $\alpha>0$, then
$$\Lambda^{p,q}_\alpha(f) \lesssim L^{p,q}_\alpha(f) + \|f\|_p.$$
\end{prop}

\begin{dem}

According to Lemma \ref{Fen3LambdaBoundedByL}, it is sufficient to check that
$$\left( \int_G \left( \frac{\|\nabla_y f\|_p e^{-c|y|^2}}{|y|^\alpha} \right)^q \frac{dy}{V(|y|)} \right)^{\frac{1}{q}}  \lesssim L^{p,q}_\alpha(f) + \|f\|_p.$$
Since we obviously have
$$\left( \int_{|y|\leq 1} \left( \frac{\|\nabla_y f\|_p e^{-c|y|^2}}{|y|^\alpha} \right)^q \frac{dy}{V(|y|)} \right)^{\frac{1}{q}} \leq L^{p,q}_\alpha(f),$$
all we need to prove is
$$T = \left( \int_{|y|\geq 1} \left( \frac{\|\nabla_y f\|_p e^{-c|y|^2}}{|y|^\alpha} \right)^q \frac{dy}{V(|y|)} \right)^{\frac{1}{q}} \lesssim \|f\|_p.$$
Indeed, $\|\nabla_y f\|_p \leq 2 \|f\|_p$ and thus
\[\begin{split}
   T & \lesssim  \|f\|_p \left( \int_{|y|\geq 1} \left( e^{-c|y|^2}\right)^q dy \right)^{\frac{1}{q}} \\
& \lesssim \|f\|_p \left( \sum_{j=0}^\infty e^{-cq4^j} V(2^{j+1}) \right)^\frac{1}{q} \\
& \lesssim \|f\|_p,
  \end{split}\]
where the last line holds because $V(r)$ have at most exponential growth. 
\end{dem}

\begin{prop} \label{Fen3LbyLambda}
Let $p,q\in [1,+\infty]$ and $\alpha \in (0,1)$. Then
$$L^{p,q}_\alpha(f) \lesssim \Lambda^{p,q}_\alpha(f) + \|f\|_p \qquad \forall f\in B^{p,q}_{\alpha}(G).$$
\end{prop}

\begin{dem}
\begin{enumerate}
 \item {\bf Decomposition of $f$:}

The first step is to decompose $f$ as
$$ f = (f-H_1 f) + H_1 f.$$

We introduce
$$f_n = -\int_{2^n}^{2^{n+1}} \frac{\dr H_t f}{\dr t} dt = -\int_{2^n}^{2^{n+1}} \Delta H_t f dt$$
and
$$c_n = \int_{2^{n-1}}^{2^{n}} \left\|\frac{\dr H_t f}{\dr t}\right\|_p dt.$$
Remark then that
$$\|f_n\|_p \leq c_{n+1}$$
and Lemma \ref{Fen3IdentityFormula} provides
$$f-H_1 f = \sum_{n=-\infty}^{-1} f_n \qquad \text{in } \mathcal S'(G).$$

\item {\bf  Estimate of $X_if_n$: }

Let  us  prove that if $n\leq -1$, one has for all $i\in \bb 1,k \bn$
\begin{equation} \label{Fen3trucque1} \|X_i f_n\|_p \lesssim 2^{-\frac n2} c_n \end{equation}
Indeed, notice first
\[\begin{split}
   f_n & = -2\int_{2^{n-1}}^{2^n} \Delta H_{2t} f dt \\
& = -2 H_{2^{n-1}}\int_{2^{n-1}}^{2^n} H_{t-2^{n-1}}\Delta H_t f dt \\
& := H_{2^{n-1}}  g_n.
  \end{split}\]
Proposition \ref{Fen3AnalycityofHt} implies then
\[\begin{split}
   \|X_i f_n\|_p & \lesssim 2^{-\frac n2} \|g_n\|_p \\
& \lesssim 2^{-\frac n2} \int_{2^{n-1}}^{2^n} \left\|H_{t-2^{n-1}}\frac{\dr H_{t} f}{\dr t}\right\|_p dt \\
& \lesssim 2^{-\frac n2} \int_{2^{n-1}}^{2^n} \left\|\frac{\dr H_{t} f}{\dr t}\right\|_p dt \\
& \qquad = 2^{-\frac n2} c_n.
  \end{split}\]

If $\varphi: \, [0,1]\to G$ is an admissible path linking $e$ to $y$  with $l(\varphi)\leq 2|y|$, 
\[\begin{split}
   \nabla_y f_n(x) & = \int_0^1 \frac{d}{ds} f_n(x\varphi(s)) ds \\
& = \int_0^1 \sum_{i=1}^k c_i(s) X_i f_n(x\varphi(s)) ds.
  \end{split}\]
Hence, \eqref{Fen3trucque1} implies
\[\begin{split}
   \|\nabla_y f_n\|_p & \leq \int_0^1 \sum_{i=1}^k |c_i(s)| \|X_i f_n(.\varphi(s))\|_p ds \\
& \qquad =   \|X_i f_n\|_p \int_0^1 \sum_{i=1}^k |c_i(s)| ds  \\
& \lesssim 2^{-\frac n2} c_n \int_0^1 \sum_{i=1}^k |c_i(s)| ds \\
& \lesssim |y|2^{-\frac n2} c_n
  \end{split}\]
where the second line is a consequence of the right-invariance of the measure  and the last one follows from the definition of $l(\varphi)$. Thus, one has
\begin{equation} \label{Fen3nablafn}
\|\nabla_y f_n\|_p \lesssim \left\{\begin{array}{ll} |y|2^{-\frac{n}{2}} c_n & \text{ if } |y|^2 < 2^n \\ c_{n+1} & \text{ if } |y|^2 \geq 2^n \end{array} \right. 
\end{equation}

\item {\bf Estimate of $L^{p,q}_\alpha(f-H_1 f)$}

As a consequence  of \eqref{Fen3nablafn}, 
\[\begin{split}
   \left[L^{p,q}_\alpha(f- H_1 f) \right]^q & = \sum_{j=-\infty}^{-1} \int_{2^j<|y|^2\leq 2^{j+1}} \left(\frac{\|\nabla_y f\|_p}{|y|^\alpha}\right) ^q \frac{dy}{V(|y|)} \\
& \lesssim \sum_{j=-\infty}^{-1} \int_{2^j<|y|^2\leq 2^{j+1}} \left(\sum_{n=-\infty}^{-1} \frac{\|\nabla_y f_n\|_p}{|y|^\alpha}\right) ^q \frac{dy}{V(|y|)} \\
& \lesssim \sum_{j=-\infty}^{-1}2^{-\frac{j\alpha}{2}q}  \left(\sum_{n=-\infty}^{j} c_{n+1} +  \sum_{n=j+1}^{-1} 2^\frac{j-n}{2} c_n\right) ^q  \\
& \lesssim \sum_{j=-\infty}^{-1}2^{jq(1-\frac{\alpha}{2})}  \left(\sum_{n=-\infty}^{-1} 2^{-\frac{\max\{j,n\}}{2}} [c_{n+1} + c_n]\right) ^q  \\
& \lesssim \sum_{n=-\infty}^{-1} \left( 2^{-\frac{n\alpha}{2}}  [c_{n+1} + c_n]\right) ^q  \\
& \lesssim \sum_{n=-\infty}^{0} \left( 2^{-\frac{n\alpha}{2}}  c_n\right)^q  \\
\end{split} \] 
 Note that the third line holds since $2^{j}\leq 1$, so that $V(2^{j+1})\lesssim V(2^j)$ and the fifth one is obtained with Lemma \ref{Fen3Bcalculus}, since $\alpha \in (0,1)$.

However
\begin{equation}\label{Fen3discretetocontinuoust}\begin{split}
   \sum_{n=-\infty}^{0} \left[ 2^{-n\frac{\alpha}{2}} c_n \right]^q & =  \sum_{n=-\infty}^{0} \left[ 2^{-n\frac{\alpha}{2}} \int_{2^{n-1}}^{2^{n}}\left\|\frac{\dr H_t f}{\dr t} \right\|_p dt \right]^q \\
& \lesssim \sum_{n=-\infty}^{0} 2^{-nq\frac{\alpha}{2}} 2^{n(q-1)} \int_{2^{n-1}}^{2^{n}}\left\|\frac{\dr H_t f}{\dr t} \right\|_p^q dt \\
& \lesssim \sum_{n=-\infty}^{0}  \int_{2^{n-1}}^{2^{n}}\left( t^{1-\frac{\alpha}{2}}\left\|\frac{\dr H_t f}{\dr t} \right\|_p\right)^q \frac{dt}{t} \\
& \qquad = \int_{0}^{1}\left( t^{1-\frac{\alpha}{2}}\left\|\frac{\dr H_t f}{\dr t} \right\|_p\right)^q \frac{dt}{t} \\
& \qquad =  \left(\Lambda^{p,q}_\alpha(f)\right)^q. 
  \end{split}\end{equation}

\item {\bf Estimate of $L^{p,q}_\alpha(H_1 f)$}

With  computations similar to those of  the second step of this proof, we find that
$$\|\nabla_y H_1 f\|_p \lesssim |y|\|f\|_p.$$
Consequently,
\[\begin{split}
    L^{p,q}_\alpha(H_1 f)  & \leq \|f\|_p \left( \int_{|y|\leq 1} |y|^{q(1-\alpha)} \frac{dy}{V(|y|)} \right)^\frac1q \\
& \leq  \|f\|_p \left( \sum_{j\leq -1} \int_{2^{j}<|y|\leq 2^{j+1}} |y|^{q(1-\alpha)} \frac{dy}{V(|y|)} \right)^\frac1q \\
& \lesssim \|f\|_p \left( \sum_{j\leq -1}  2^{qj(1-\alpha)} \right)^\frac1q \\
& \lesssim \|f\|_p
  \end{split}\]
where the third line is a consequence of the local doubling property.
\end{enumerate}
\end{dem}

\begin{theo}
Let $G$ be a unimodular Lie group and $\alpha \in (0,1)$, then we have the following Leibniz rule.

If $p_1,p_2,p_3,p_4,p,q \in [1,+\infty]$ are such that
$$\frac{1}{p_1} + \frac{1}{p_2} = \frac{1}{p_3} + \frac{1}{p_4} = \frac{1}{p}$$
then for all $f\in B^{p_1,q}_\alpha(G) \cap L^{p_3}(G)$ and all $g \in B^{p_4,q}_\alpha(G) \cap L^{p_2}(G)$, one has
$$\|fg\|_{B^{p,q}_\alpha} \lesssim \|f\|_{B^{p_1,q}_\alpha} \|g\|_{L^{p_2}} + \|f\|_{L^{p_3}} \|g\|_{B^{p_4,q}_\alpha}.$$
\end{theo}

\begin{dem}
Check that
$$\nabla_y (f\cdot g)(x) = g(xy)\cdot \nabla_y f(x) + f(x)\cdot \nabla_y g(x).$$
Thus, with Hölder inequality,
\[\begin{split}
   \|fg\|_{B^{p,q}_\alpha} & \simeq \|f\cdot g\|_p + L^{p,q}_\alpha(f\cdot g) \\
& \lesssim \|f\|_{p_1} \|g\|_{p_2} + L^{p_1,q}_\alpha(f)\cdot \|g\|_{L^{p_2}} + \|f\|_{L^{p_3}} \cdot L^{p_4,q}_\alpha(g) \\
& \lesssim \|f\|_{B^{p_1,q}_\alpha} \|g\|_{p_2} + \|f\|_{L^{p_3}} \|g\|_{B^{p_4,q}_\alpha}.
  \end{split}\]
\end{dem}

\subsection{ Characterization  by induction - Theorem \ref{Fen3Main4}}

\begin{prop} \label{Fen3PropInduction}
Let $p,q\in [1,+\infty]$ and $\alpha > -1$. Let $m > \frac{\alpha}{2}$. One has for all $i\in \bb 1,k \bn$
$$\Lambda^{p,q}_\alpha(X_i f)  \lesssim \Lambda^{p,q}_{\alpha + 1} f + \|f\|_p = \|f\|_{B^{p,q}_{\alpha+1}}.$$
\end{prop}

\begin{dem}
The scheme of the proof is similar to Proposition \ref{Fen3DiscreteEquivalenceC}. 

\begin{enumerate}
 \item {\bf Decomposition of $f$:}

Let $M$ be  an integer with $M>\frac{\alpha + 1}{2}$. 
We decompose $f$ as  in  Lemma \ref{Fen3IdentityFormula}:
$$ f = \frac{1}{(M-1)!} \int_0^{1} (t\Delta)^{M} H_t f \frac{dt}{t} +  \sum_{k=0}^{M-1} \frac{1}{k!} \Delta^k H_1 f$$ 
and we introduce
$$f_n = -\int_{2^n}^{2^{n+1}} (t\Delta)^M  H_tf \frac{dt}{t}$$
and
$$c_n = \int_{2^{n-1}}^{2^{n}} t^M\left\|\Delta^M H_t f\right\|_p \frac{dt}{t}.$$
Remark then that
$$\|f_n\|_p \leq c_{n+1}$$
and
$$f = \frac{1}{(M-1)!}\sum_{n=-\infty}^{-1} f_n + \sum_{k=0}^{M-1} \frac{1}{k!} \Delta^k H_1 f.$$

\item {\bf  A first estimate of $\Delta^mH_tX_if_n$: }

Let us prove that if $n\leq -1$, one has for all $i\in \bb 1,k \bn$
\begin{equation} \label{Fen3trucque2} \|\Delta^m X_i f_n\|_p \lesssim 2^{-n(m+\frac 12)} c_n. \end{equation}
Indeed, notice first
\[\begin{split}
   f_n & =  -2^M  \int_{2^{n-1}}^{2^n}  (t\Delta)^M H_{2t} f \frac{dt}t \\
& = -2^M  H_{2^{n-1}}\int_{2^{n-1}}^{2^n} H_{t-2^{n-1}} (t\Delta)^M H_{t} f \frac{dt}t \\
& := H_{2^{n-1}}  g_n.
  \end{split}\]
Thus, since $\Delta = -\sum_{i=1}^k X_i^2$ can be written as a polynomial in the $X_i$'s, we obtain with the upper estimate of the heat kernel (Proposition \ref{Fen3htestimates}),
\[\begin{split}
   \|\Delta^m X_i f_n\|_p & \lesssim \left( \int_G \left|\int_G |\Delta^M X_i h_{2^{n-1}}(z^{-1}x)| g_n(z) dz\right|^p dx \right)^\frac1p  \\
& \lesssim  \frac{2^{-n(m+\frac{1}{2})}}{V(2^\frac{n}{2})} \left( \int_G \left|\int_G \exp\left(-c \frac{|z^{-1}x|^2}{2^n} \right) g_n(z)dz\right|^p dx \right)^\frac1p \\
& \lesssim 2^{-n(m+\frac{1}{2})}  \|g_n\|_p \\
& \lesssim 2^{-n(m+\frac{1}{2})} c_n
  \end{split}\]
where  the second line is due to the fact that $V(2^{\frac n2})\lesssim V(2^{\frac{n-1}2})$ and the last two  lines are obtained  by an argument analogous to the one for  \eqref{Fen3trucque1}.

$$
\begin{array}{lll}
\left\Vert \Delta^mH_tX_if_n\right\Vert_p & \lesssim &\frac 1{t^m} \left\Vert X_if_n\right\Vert_p\\
& \lesssim & \frac 1{t^m} \left\Vert X_iH_{2^{n-1}}g_n\right\Vert_p\\
\end{array}
$$

 As a consequence, one has for all $t\in (0,1]$,
$$\left\|\Delta^m H_t X_i f_n \right\|_p  = \|H_t \Delta^m X_i f_n\|_p \lesssim 2^{-n(m+\frac{1}{2})} c_n,$$
since $H_t$ is uniformly bounded.

\item {\bf A second  estimate of $\Delta^mH_tX_if_n$: }

Let us prove that for all $f \in L^p(G)$ and for all $i\in \bb 1,k \bn$, one has
\begin{equation} \label{Fen3trucque3} \left\|\Delta^m H_t X_i f\right\|_p \lesssim t^{-m-\frac{1}{2}} \|f\|_p. \end{equation}
First, notice that

\[\begin{split}
   \Delta^m H_t X_i f(x) & = \int_G \frac{\partial^m}{\partial t^m}h_t(y) (X_i f)(xy) dy \\
& = \int_G \frac{\partial^m}{\partial t^m}h_t(y) [X_i f(x.)](y) dy \\
& = -\int_G X_i \frac{\partial^m}{\partial t^m}h_t(y) f(xy) dy \\
& = -\int_G X_i \Delta^m h_t(x^{-1}y) f(y) dy. 
  \end{split}\]
Then, using the estimates on the heat kernel (Proposition \ref{Fen3htestimates})and the fact that $\Delta = -\sum X_i^2$, we obtain
\[\begin{split}
   \left\|\Delta^m  H_t X_i f\right\|_p & \lesssim \left(\int_G \left| \frac{t^{-m-\frac12}}{V(\sqrt t)} \int_G \exp\left(-c\frac{|x^{-1}y|^2}{t} \right)|f(y)| dy \right|^p dx \right)^\frac1p \\
& : = t^{-m-\frac12} \left(\int_G \left| \int_G K(x,y)|f(y)| dy \right|^p dx \right)^\frac1p
  \end{split}\]
with $K(x,y) = \frac{1}{V(\sqrt t)} \exp\left(-c\frac{|x^{-1}y|^2}{t} \right)$. Proposition \ref{Fen3convolutionbyht} yields the estimate \eqref{Fen3trucque3}.

\item {\bf Estimate of $\Lambda^{p,q}_\alpha(\sum f_n)$}

The two previous steps imply
$$\left\|\Delta^m H_t X_i f_n \right\|_p \lesssim \left\{\begin{array}{ll} 2^{-n(m+\frac{1}{2})} c_n & \text{ if } t<2^n \\ t^{-m-\frac{1}{2}}c_{n+1} & \text{ if } t \geq 2^n \end{array} \right. .$$
As a consequence,
\[\begin{split}
   \int_0^1 & \left( t^{m-\frac{\alpha}{2}} \left\|\Delta^m H_t X_i \sum_{n=-\infty}^{-1}f_n \right\|_p \right)^q   \frac{dt}{t} \\
& \qquad \lesssim \sum_{j=-\infty}^{-1} \int_{2^j<t\leq 2^{j+1}} \left( t^{m-\frac{\alpha}{2}} \sum_{n=-\infty}^{-1} \left\|\Delta^m H_t X_if_n\right\|_p \right)^q \frac{dt}{t} \\
& \qquad \lesssim \sum_{j=-\infty}^{-1} \left(2^{j(m-\frac{\alpha}{2})}\sum_{n=-\infty}^{j} 2^{-j(m+\frac12)}c_{n+1} +  \sum_{n=j+1}^{-1} 2^{-n(m+\frac12)} c_n\right) ^q  \\
& \qquad \lesssim \sum_{j=-\infty}^{-1} \left(2^{j(m-\frac{\alpha}{2})}\sum_{n=-\infty}^{-1} 2^{-\max\{j,n\}(m+\frac12)} [c_n + c_{n+1}]\right) ^q  \\
& \qquad \lesssim \sum_{n=-\infty}^{-1} \left[2^{-n\frac{\alpha + 1}{2}} [c_n+c_{n+1}]\right]^q \\
& \qquad \lesssim \sum_{n=-\infty}^0 \left[2^{-n\frac{\alpha + 1}{2}} c_n\right]^q 
\end{split} \] 
where we used Lemma \ref{Fen3Bcalculus} for the fourth estimate, relevant since $-1<\frac{\alpha}{2}<m$ by assumption.
We get then the domination
\begin{equation} \int_0^1 \left( t^{m-\frac{\alpha}{2}} \left\|\Delta^m H_t X_i \sum_n f_n \right\|_p \right)^q \frac{dt}{t} \lesssim \sum_{n=-\infty}^{0} \left[ 2^{-n\frac{\alpha + 1}{2}} c_n \right]^q.\end{equation}
However  computations analogous to those leading to \eqref{Fen3discretetocontinuoust}  prove that
\begin{equation}\begin{split}
   \sum_{n=-\infty}^{0} \left[ 2^{-n\frac{\alpha+1}{2}} c_n \right]^q &  \lesssim  \int_{0}^{1}\left( t^{M-\frac{\alpha+1}{2}}\left\|\Delta^M f\right\|_p\right)^q \frac{dt}{t} \\
& \lesssim  \left(\Lambda^{p,q}_{\alpha+1} f\right)^q .
  \end{split}\end{equation}

\item {\bf Estimate of the remaining term.}

Recall that 
$$f = \frac{1}{(M-1)!}\sum_{n=-\infty}^{-1} f_n + \sum_{k=0}^{M-1} \frac{1}{k!} \Delta^k H_1 f.$$
We  already  estimated $\Lambda^{p,q}_\alpha(\sum f_n)$.  What remains to be estimated is  $\Lambda^{p,q}_\alpha(\sum \frac{1}{k!} \Delta^k H_1 X_if)$.

Proposition \ref{Fen3AnalycityofHt}   provides as well 
$$\left\|\Delta^m H_t X_i \Delta^k H_1 f\right\|_p \leq \left\|\Delta^m X_i \Delta^k H_1 f \right\|_p  \lesssim \|f\|_p.$$

As a consequence, we get,
\[\begin{split} \int_0^1 \left( t^{m-\frac{\alpha}{2}} \left\|\Delta^m H_t X_i \sum_{k=0}^{M-1} \frac{1}{k!} \Delta^k H_1 f \right\|_p \right)^q \frac{dt}{t}
&  \lesssim  \left(\|f\|_p \int_0^1  t^{q(m-\frac{\alpha}{2})}  \frac{dt}{t}\right)^q \\
& \lesssim \|f\|_p^q.
  \end{split}\]
\end{enumerate}
\end{dem}

\begin{cor} Let $p,q\in [1,+\infty]$ and $\alpha >0$.

$$\|f\|_{B^{p,q}_{\alpha+1}} \simeq \|f\|_{L^p} + \sum_{i=1}^k \|X_i f\|_{B^{p,q}_\alpha}.$$
\end{cor}

\begin{dem}
The main work was done in the previous proposition. Indeed, notice that Proposition \ref{Fen3PropInduction} implies
\[\begin{split}
\Lambda^{p,q}_{\alpha+1} f & = \Lambda^{p,q}_{\alpha-1} (\Delta f) \\
& \leq \sum_{i=1}^k \Lambda^{p,q}_{\alpha-1} X_i(X_if) \\
& \lesssim \sum_{i=1}^k \|X_i f\|_{B^{p,q}_\alpha},
  \end{split}\]
which provides the domination of the first term by the second one.

The  converse inequality splits into  two parts. 
The first one is the domination of $\Lambda^{p,q}_\alpha (X_i f)$ by $\|f\|_{B^{p,q}_{\alpha+1}}$, which is  an  immediate application of Proposition \ref{Fen3PropInduction}.
The second one is the domination of $\|X_i f\|_p$. 
But recall that Theorem \ref{FEN3MAIN1} states that we can replace $\|X_i f\|_p$ by $\|H_\frac{1}{2} X_i f\|_p$  in the Besov norm, and  \eqref{Fen3trucque3} provides that
$$\|H_\frac{1}{2} X_i f\|_p \lesssim \|f\|_p \leq \|f\|_{B^{p,q}_\alpha}.$$ 
\end{dem}

\bibliographystyle{plain}
\bibliography{Biblio}

\end{document}